\documentclass[11pt,reqno]{amsart}

\usepackage{euscript}
\usepackage{pb-diagram}
\usepackage{lamsarrow}
\usepackage{pb-lams}
\usepackage{epic}
\usepackage{epsfig}
\usepackage{amssymb}

\theoremstyle{plain}

\theoremstyle{definition}

\newcommand\bC{{\mathbb C}}

\newcommand{\bQ}{{{\mathbb Q}}}

\newcommand{\Z}{\EuScript{Z}}

\newcommand\aug{\mathfrak{aug}}
\newcommand{\C}{\EuScript{C}}
\newcommand{\p}{\EuScript{P}}

\newcommand{\ra}{\rightarrow}

\newcommand{\lan}{\langle}
\newcommand{\ran}{\rangle}

\def\inpr{\mathbin{\hbox to 6pt{\vrule height0.4pt width5pt depth0pt \kern-.4pt \vrule height6pt width0.4pt depth0pt\hss}}}

\newcommand{\si}{{\sigma}}

\newcommand{\eps}{\epsilon}
\newcommand{\ssw}{{\bf sw}}
\newcommand{\lag}{\langle}
\newcommand{\rag}{\rangle}
\newcommand{\no}{\natural}

\newcommand{\et}{\EuScript{T}}
\newcommand{\M}{\EuScript{M}}

\newcommand{\bms}{\mbox{\boldmath$s$}}

\newcommand{\sign}{\mbox{sign}}

\def\mmod{\mbox{mod}}

\def\C{\mathbb C}
\def\Q{\mathbb Q}

\def\Z{\mathbb Z}
\def\N{\mathbb N}
\begin{document}

\title{Seiberg-Witten invariants and surface singularities III\\
(splicings and cyclic covers)}

\author{Andr\'as N\'emethi}
\address{Department of Mathematics\\Ohio State University\\Columbus, OH 43210}
\email{nemethi@math.ohio-state.edu}
\urladdr{http://www.math.ohio-state.edu/ \textasciitilde nemethi/}
\author{Liviu I. Nicolaescu}
\address{University of Notre Dame\\Notre Dame, IN 46556}
\email{nicolaescu.1@nd.edu} \urladdr{http://www.nd.edu/
\textasciitilde lnicolae/}
\thanks{The first author is partially supported by NSF grant DMS-0088950;
the second author is partially supported by NSF grant
DMS-0071820.}

%\date{December, 2001}

%\date{Version 1, February 2001}

\begin{abstract} We verify the conjecture formulated in \cite{NN} for
suspension  singularities of type $g(x,y,z)=
f(x,y)+z^n$, where $f$ is an irreducible
plane curve singularity. More precisely, we prove that the modified
Seiberg-Witten invariant of the link $M$ of $g$, associated with the
canonical $spin^c$ structure, equals $-\sigma(F)/8$, where $\sigma(F)$ is the
signature of the Milnor fiber of $g$.  In order to do this, we prove
general
splicing formulae for the Casson-Walker invariant and for the sign refined
Reidemeister-Turaev torsion (in particular, for the modified
Seiberg-Witten invariant too).
These provide results for some cyclic covers as well.
As a by-product, we compute all the relevant invariants of $M$
in terms of the Newton pairs of $f$ and the integer $n$.
\end{abstract}

\maketitle
\pagestyle{myheadings}
\markboth{{\normalsize Andr\'as N\'emethi and Liviu I. Nicolaescu}}{
{\normalsize Seiberg-Witten invariants and surface singularities III}}

\section{Introduction}

The present article is a natural continuation of \cite{NN} and \cite{NNII},
and it is closely related to \cite{nemsignat,NDedI,NDedII}.
In \cite{NN}, the authors formulated  a very
general conjecture  which connects  the topological and the
analytical invariants  of a complex normal surface singularity
whose link is  a rational homology sphere.

Even if we restrict ourselves to the case of
hypersurface singularities, the conjecture is still highly non-trivial.
The ``simplified'' version for  this  case reads as follows.

\subsection{Conjecture.}\label{i1} \cite{NN} \ {\em Let $g:(\C^3,0)\to
(\C,0)$ be a complex  analytic  germ which defines an isolated hypersurface
singularity. Assume that its link $M$ is a rational homology sphere.
Denote by $\ssw^0_M(\sigma_{can})$  the modified Seiberg-Witten invariant
of $M$ associated with the canonical $spin^c$ structure $\sigma_{can}$
(cf. \ref{sw}). Moreover, let $\sigma(F)$ be the signature of the Milnor fiber
$F$ of $g$. Then }
\begin{equation*}
-\ssw^0_M(\sigma_{can})=\sigma(F)/8.
\tag{1}
\end{equation*}

The goal of the present paper is to verify this conjecture for suspension
hypersurface singularities. More precisely, in \ref{conj} we prove the
following.

\subsection{Theorem.}\label{i2} {\em Let $f:(\C^2,0)\to (\C,0)$ be an
irreducible plane curve singularity. Fix an arbitrary positive integer $n$
such that the link $M$ of the suspension singularity $g(x,y,z):=f(x,y)+z^n$
is a rational homology sphere (cf. \ref{pr}(3)). Then \ref{i1}(1) holds.}

\vspace{2mm}

The numerical identity \ref{i1}(1) covers a very deep qualitative
analytic-rigidity phenomenon.

From topological point of view,
any normal two-dimensional analytic  singularity $(X,0)$ is
completely characterized by its link $M$, which is an oriented
3-manifold. Moreover, by a result of Neumann \cite{NP}, any
decorated plumbing (or resolution)
graph of $(X,0)$ carries the same information
as $M$. A property of $(X,0)$ will be called \emph{topological} if
it can be determined from  any of these graphs.

A very intriguing issue, which has generated intense research efforts, is  the possibility of expressing the
analytic invariants of $(X,0)$ (like the geometric genus $p_g$,
multiplicity, etc.) or the smoothing invariants (if they exist, like
the signature $\sigma(F)$ or the topological Euler-characteristic $\chi(F)$
of  the Milnor fiber $F$) in terms of the topology of $M$.

\subsection{A short historical survey.}\label{i3} \
M. Artin  proved in \cite{Artin,Artin2}
that the rational singularities (i.e. the vanishing of $p_g$) can be
characterized  completely from the resolution (or plumbing) graph.
In \cite{Lauferme}, H. Laufer   extended  Artin's results to minimally elliptic
singularities, showing that Gorenstein singularities with $p_g=1$ can be
characterized topologically.
Additionally, he noticed that the program breaks for
more complicated singularities (see also the comments in \cite{NN} and
\cite{Neminv}).
On the other hand, the first author noticed in \cite{Neminv} that
Laufer's counterexamples do not signal  the end of the program. He
conjectured  that if we restrict ourselves to the case of those
Gorenstein singularities whose links are rational homology spheres
then  $p_g$ is topological. This was carried out explicitly for elliptic
singularities in \cite{Neminv} (partially based on some results of
S. S.-T. Yau, cf.  e.g.  with \cite{Yau}).

For Gorenstein singularities which have smoothing
(with Milnor fiber $F$), the topological invariance of $p_g$
can be  reformulated in terms of $\si(F)$ and/or $\chi(F)$. Indeed,
via some results of Laufer, Durfee, Wahl  and Steenbrink, any of $p_g$,
$\si(F)$ and $\chi(F)$ determines  the remaining two  modulo
$K^2+\#{\mathcal V}$ (for the precise identities see e.g.  \cite{LW}
or \cite{NN}).
Here, $K^2+\#{\mathcal V}$ is  defined as follows.
For a given resolution, one  takes
the canonical divisor $K$, and the number $\# {\mathcal V}$
of irreducible components of the exceptional divisor of the resolution.
Then $K^2+\#{\mathcal V}$ can be deduced from the resolution graph,
and is independent of the choice of the
graph. In particular, it is an invariant of the link $M$.

For example, the identity  which connects $p_g$ and $\sigma(F)$ is
\begin{equation*}
8p_g+\sigma(F)+K^2+\#{\mathcal V}=0.
\tag{2}
\end{equation*}
This  connects the above facts about $p_g$
with the following list of results about $\si(F)$.

Fintushel and Stern proved in \cite{FS}  that
for a hypersurface Brieskorn  singularity whose link is an {\em
integral} homology sphere, the Casson invariant $\lambda(M)$  of
the link $M$ equals $\si(F)/8$.  This fact was generalized
by Neumann  and Wahl in \cite{NW}. They proved the same statement
for all Brieskorn-Hamm complete intersections and suspensions of
irreducible  plane
curve singularities (with the same assumption about the link).
Moreover, they conjectured  the validity
of the formula for any isolated complete intersection  singularity
whose link is an {\em integral}  homology sphere.

In \cite{NN} the authors extended  the above conjecture for
smoothing of Gorenstein singularities with {\em rational} homology sphere link.
Here   the Casson invariant $\lambda(M)$ is replaced by
a certain Seiberg-Witten invariant $\ssw^0_M(\sigma_{can})$
of the link associated with the {\em canonical} $spin^c$ structure of $M$.
$\ssw^0_M(\sigma_{can})$ is the difference of a certain
Reidemeister-Turaev sign-refined torsion invariant and the
Casson-Walker invariant, for details see \ref{sw}. If $M$ is an integral
homology sphere, then the torsion invariant vanishes.

In fact, the conjecture in \cite{NN} is more general.
A part of it says that for any $\Q$-Gorenstein singularity whose link
is a rational homology sphere, one has
\begin{equation*}
8p_g-8\ssw^0_M(\sigma_{can})+K^2+\#{\mathcal V}=0.
\tag{3}
\end{equation*}
Notice that (in the presence of a smoothing and of the Gorenstein property,
e.g. for any hypersurface singularity)
(3) via (2) is exactly (1). The identity (3)
was verified in \cite{NN} for cyclic quotient singularities, Brieskorn-Hamm
complete intersections and some rational and minimally elliptic singularities.
\cite{NNII} contains the case when $(X,0)$ has a good $\C^*$-action.

Finally we mention, that recently Neumann and Wahl initiated in \cite{NW2} an
all-embracing program about those $\Q$-Gorenstein singularities whose
link is rational homology sphere. They conjecture that the universal abelian
cover of such a singularity is an isolated complete intersection, and
its Milnor fiber  can be recovered
from $M$ together with the action of $H_1(M)$ (modulo some equisingular and
equivariant deformation).  In particular, this implies the topological
invariance  of $p_g$ as well. We hope that our
efforts in the direction of the conjecture formulated in \cite{NN}
 will contribute to the accomplishment of this  program as well.

\vspace{2mm}

On the other hand, the theory of suspension hypersurface
singularities also has its own long history.  This class  (together
with the weighted-homogeneous singularities) serve as an important
``testing and exemplifying'' family for various properties and conjectures.
For more information, the reader is invited to check \cite{NDedI,NDedII}
and the survey paper \cite{nemsignat}, and the references listed in these
articles.

\subsection{A few words about the proof.}\label{i4} \ If $M$ is the
link of $g=f+z^n$ (as in \ref{i2}), then $M$ has a natural splice
decomposition  into Seifert varieties of type $\Sigma(p,a,m)$. Moreover,
in \cite{NDedI}(3.2) the first author established an additivity formula
for $\sigma(F)$ compatible with the geometry of this decomposition.
On the other hand, for any Brieskorn singularity $(x,y,z)\mapsto
x^p+y^a+z^m$ (whose link is
$\Sigma(p,a,m)$) the conjectured identity (1) is valid by \cite{NN,NNII}.
Hence it was natural to carry out the proof of \ref{i2} by proving
an  additivity result for
$\ssw^0_M(\sigma_{can})$ with respect to the splice decomposition of $M$
into Seifert varieties.

This additivity result is proved in \ref{swadd} (as an outcome of all the
preparatory results  of the previous sections) but its proof contains
some surprising steps.

Our original plan was the following. First, we identify the splicing
data of $M$. Then, for such  splicing data, we establish
splicing formulas for the Casson-Walker invariant
and for the Reidemeister-Turaev
sign-refined torsion  with the hope that we can do this
in the world of topology without going back for some extra restrictions
to the world of singularities. This program  for the first invariant
was   straightforward, thanks to the results of Fujita \cite{fujita}
and Lescop \cite{Lescop} (see section 3).
But when dealing with the torsion we encountered some serious difficulties
(and finally we had to return back to singularities for some additional
properties).

The torsion-computations  require the explicit description of
the supports of all the relevant characters, and then the computation of  some
sophisticated Fourier-Dedekind sums.  The computation turned out to be
feasible because these  sums are not arbitrary. They have two very subtle
special features  which  follow from various properties of irreducible plane
curve singularities. The first one is a numerical inequality
(\ref{a1}(6))
(measuring some strong algebraic rigidity). The second (new)
property is the alternating
property of their Alexander polynomial (\ref{a2}).

In section 4
we establish different splicing formulas for $\ssw_M^0(\sigma_{can})$, and we
show the limits of a possible additivity. We introduce even a new invariant
${\mathcal D}$ which measures the non-additivity property of
$\ssw_M^0(\sigma_{can})$ with respect to (some) splicing (see e.g. \ref{g9})
or (some) cyclic covers (see \ref{swc}). This invariant vanishes in the
presence of the alternating property of the involved Alexander polynomial.

This shows clearly (and rather surprisingly) that the behavior of
$\ssw_M^0(\sigma_{can})$
with respect to splicing and cyclic covers (constructions topological in
nature) definitely prefers some special
algebraic situations. For more comments,
see \ref{g12}, \ref{a4} and \ref{rem}.

Section 5 contains the needed  results for irreducible plane curve
singularities  and the Algebraic Lemma used in the summation of the
Fourier-Dedekind sums mentioned above.

In section 6 we provide a list of properties
of the link of $f+z^n$. Here basically we use almost all the partial results
proved in the previous sections. Most of the formulae are formulated as
inductive  identities with respect to the number of Newton pairs of $f$.

\subsection{Notations.} \
All the homology groups  with unspecified coefficients are
defined over the integers. In section 3,  $\bms(\cdot,\cdot)$ denotes
Dedekind sums, defined  by the same convention as in \cite{fujita, Lescop}
or \cite{NN}. For definitions and detailed discussions about the involved
invariants, see \cite{NN}. For different
 properties of hypersurface singularities, the reader may consult \cite{AGV}
as well.

\section{Preliminaries and notations}

\subsection{Oriented knots in rational homology spheres.}\label{s0} Let $M$ be an
oriented 3-manifold which is a {\em rational homology sphere}.
Fix an oriented  knot $K\subset M$, denote by $T(K)$ a small tubular
neighborhood  of $K$ in $M$, and let $\partial T(K)$ be its
oriented boundary with its natural orientation.
The natural {\em oriented meridian} of $K$, situated in
$\partial T(K)$,  is denoted be $m$.
We fix an oriented {\em parallel} $\ell$ in $\partial T(K)$ (i.e.
$\ell\sim K$ in $H_1(T(K))$). If $\lag,\rag$
denotes the intersection form
in $H_1(\partial T(K))$, then $\lag m,\ell\rag=1$
 (cf. e.g. Lescop's book \cite{Lescop}, page 104;
we will use the same notations $m$ and $l$ for some geometric realizations
of the meridian and parallel as primitive simple curves,
respectively for their homology classes
in $H_1(\partial T(K))$).

Obviously, the choice of $\ell$ is not unique.
In all our applications, $\ell$ will be characterized by some
precise additional geometric construction.

Assume that the order of the homology class of $K$ in $H_1(M)$ is
$o>0$. Consider an oriented surface $F_{o K}$ with boundary
$oK$, and take the intersection $\lambda:=
F_{o K}\cap \partial T(K)$.  $\lambda$ is called the {\em
longitude} of $K$.
The homology class of $\lambda$
in $H_1(\partial T(K))$ can be represented as
$\lambda=o\ell+km$ for some integer $k$. Set
gcd$(o,|k|)=\delta>0$. Then $\lambda$ can be represented in
$\partial T(K)$ as $\delta$ primitive  torus curves  of type
$(o/\delta,k/\delta)$ with respect to $\ell$ and   $m$.

%Sometimes it is convenient (see e.g. \cite{fujita}, and \ref{s1b} and \ref{s2}
%below)  to introduce a basis
%$(x,y)$ in $H_1(\partial T(K))$  so that  $\delta y=\lambda$ and
%$\lag x,y\rag=1$. Obviously, $y=(o\ell +km)/\delta$.
%Assume that
%$x$ has the form $u\ell+vm$ for some integers $u$ and $v$.
%Then $\lag x,y\rag=1$ reads as
%\begin{equation*}
%vo-uk=\delta.\tag{1}\end{equation*}

\subsection{Dehn fillings.}\label{s0b}  Let $T(K)^\circ$ be the interior of
$T(K)$. For any homology class $a\in H_1(\partial T(K))$,
which can be represented by a primitive simple closed curve in
$\partial T(K)$, one defines the {\em Dehn filling} of
$M\setminus T(K)^\circ$ along $a$ by
$$(M\setminus T(K)^\circ )(a)=M\setminus T(K)^\circ\coprod_f S^1\times D^2,$$
where $f:\partial (S^1\times D^2)\to \partial T(K)$ is a diffeomorphism
which sends $\{*\}\times \partial D^2$ to a curve representing $a$.

\subsection{Linking numbers.}\label{s0c}
Consider two
oriented knots  $K, \ L\subset  M$ with $K\cap L=\varnothing$. Fix a
Seifert surface $F_{oK}$ of $oK$ (cf. \ref{s0})
and define  the  linking number
$Lk_M(K,L)\in \Q$ by the ``rational'' intersection $(F_{oK}\cdot L)/o$.
In fact, $Lk_M(K,\cdot):H_1(M\setminus K,\Q)\to \Q$ is a well-defined
homeomorphism and $Lk_M(K,L)=Lk_M(L,K)$.
For any  oriented knot  $L$ on $\partial T(K)$  one has
(see e.g. \cite[6.2.B]{Lescop}):
\begin{equation*}
Lk_M(L,K )=\lan L,\lambda\ran/o.
\tag{1}\end{equation*}
For any oriented knot  $K\subset M$  one has the obvious exact sequence
\begin{equation*}
0\to \Z\stackrel{\alpha}{\longrightarrow}
H_1(M\setminus T(K))\stackrel{j}{\longrightarrow} H_1(M)\to 0,
\tag{2}
\end{equation*}
where $\alpha(1_\Z)=m$. If $K\subset M$ is homologically trivial then this
sequence splits. Indeed, let $\phi$ be the restriction of
$Lk_M(K,\cdot)$  to $H_1(M\setminus K)=H_1(M\setminus T(K))$. Then $\phi$
has integer values and  $\phi\circ \alpha=1_\Z$. This provides automatically
a morphism $s:H_1(M)\to H_1(M\setminus T(K))$ such that
$j\circ s=1$ and $\alpha\circ \phi+s\circ j=1$; in particular
with  $\phi \circ s=0$ too. In fact, $s(H_1(M))=Tors\, H_1(M\setminus T(K))$.
Moreover, under the same assumption $o=1$, one has the isomorphisms
\begin{equation*}
H_2(M,K)\stackrel{\partial }{\to } H_1(K)=\Z\ \ \mbox{and}\ \
H_1(M)\to H_1(M,K).
\tag{3}
\end{equation*}

Sometimes, in order to simplify the notations, we write $H$ for the group
$H_1(M)$.

The finite group $H$ carries a natural symmetric bilinear form
$b_M:H^{\otimes 2}\to \Q/\Z$ defined by $b_M([K],[L])=Lk_M(K,L) \
(\mmod\ \Z) $, where $L$ and $K$ are two representatives with $K\cap L=
\varnothing$. If $\hat{H}$ denotes the Pontryagin dual $Hom(H,S^1)$
of $H$, then $[K]\mapsto b_M([K],\cdot)$ is an isomorphism
$H\to \hat{H}$.

\subsection{$(M,K)$ represented by plumbing.}\label{s0d}
The main application of the present article involves algebraic links $(M,K)$
which can be represented by plumbing. We recall the notations briefly
(for more details, see e.g. \cite{NN} or \cite{NeR}).

We will denote by $\Gamma(M,K)$ the plumbing graph of a link $K\subset M$.
The vertices
$v\in {\mathcal V}$ are decorated by the Euler numbers $e_v$ (of the
$S^1$-bundles over $E_v\approx S^2$ used in the plumbing construction).
The components of the link $K$ are represented by arrows
in $\Gamma(M,K)$:
if an arrow  is attached to the vertex $v$ then the corresponding
component of $K$ is a fixed fiber of the $S^1$-bundle over $E_v$.
(We think about an arrow as an arrowhead connected to $v$ by
an edge.)
If we delete the arrows then we obtain a plumbing graph $\Gamma(M)$ of $M$.
Let $\delta_v$ (resp. $\bar{\delta}_v$) be the degree (i.e. the number of
incident edges) of the vertex $v$ in $\Gamma(M)$ (resp. in $\Gamma(M,K)$).
Evidently $\bar{\delta}_v-\delta_v$ is exactly the number of arrows
supported by the vertex $v$.

Since $M$ is a rational homology sphere, $\Gamma$ is a tree.

Let $\{I_{uv}\}_{u,v\in {\mathcal V}}$ be the intersection matrix
associated with $\Gamma$; i.e. $I_{uu}=e_v$,  and for $u\not=v$ the entry
$I_{uv}=$ 0 or 1 depending that $u$ and $v$ are connected or not
in  $\Gamma$. Since $M$ is rational homology sphere,
 $I$ is non-degenerate. In fact
\begin{equation*}
|\det(I)|=|H|.
\tag{4}
\end{equation*}

The generic fiber of the $S^1$-bundle over $E_v$
is denoted by $g_v$, and we use the
same notation for its homology class in $H_1(M)$ as well.
By the above discussion, if $u\not=v$ then $Lk_M(g_u,g_v)$ is well-defined.
If $u=v$ then we  write $Lk_M(g_u,g_u) $ for $Lk_M(g_u,g_u')$ where
$g_u$ and $g_u'$ are two different fibers of the $S^1$-bundle over $E_v$.

For any fixed vertex $u\in {\mathcal V}$, we denote
by $\vec{b}(u)$ the column vector with entries =1 on the place $u$ and
zero otherwise. We define the column vector
$\vec{w}(u)$ (associated with the knot $g_{u}\subset  M$ and its order $o(u)$)
as the solution of the (non-degenerate) linear
system  $I\cdot \vec{w}(u)=-o(u)\vec{b}(u)$. The entries $\{w_{v}(u)\}_{
v\in {\mathcal V}}$ are called the weights associated with $g_{u}$.
Then the inverse matrix $I^{-1}$ of $I$, the set of weights $\{w_{v}(u)\}_{
v\in {\mathcal V}}$, and the linking pairing $Lk_M$  satisfy:
\begin{equation*}
-I^{-1}_{uv}=\frac{w_v(u)}{o(u)}=Lk_M(g_u,g_v) \ \ \
\mbox{for any $u$ and $v\in {\mathcal V}$}.
\tag{5}
\end{equation*}
In particular, $w_v(u)\in \Z$ for any $u$ and $v$.
In fact, if $I$ is negative definite, then the integers $w_v(u)$ are all
positive.

\subsection{$(M,K)$ represented by splice diagram.}\label{spl} If $M$ is an
{\em integral homology sphere}, and $(M,K)$ has a plumbing representation,
then there is an equivalent graph-codification  of $(M,K)$ in terms
of the {\em splice} (or {\em  Eisenbud-Neumann}) {\em diagram}, for details
see \cite{EN}.

The splice diagram preserves the ``shape'' of the plumbing graph (e.g. there
is a one-to-one correspondence between those vertices $v$ with $\delta_v\not=2$
 of the splice, respectively of the plumbing  graphs), but in the splice
diagram one collapses into an edge each  string of the plumbing graph.
Moreover, the decorations are also different. In the splice
diagram, each vertex has a sign $\epsilon =\pm 1$, which in all our cases
will be $\epsilon=+1$, hence we omit them. Moreover, if an end of an edge
is attached to a vertex  $v$ with $\delta_v\geq 3$, then it has a positive
integer as its decoration. The arrows have the same significance.

One of the big advantages of the of the splice diagram is that it codifies
in an ideal way the splicing decomposition of $M$ into Seifert pieces.
(In fact, the numerical decorations are exactly the Seifert invariant of
the corresponding Seifert splice-components.)

Therefore, in some cases it is much easier and more suggestive to use them.
(Nevertheless, we will use them only in those cases when we really want to
emphasize this principle, e.g. in the proof of  \ref{swc}, or when it
is incomparably easier to describe a construction with them, e.g. in
\ref{dnotzero}.)
The reader is invited to consult the book of Eisenbud and Neumann \cite{EN}
for the needed properties. The criterion which
guarantees that $(M,K)$ is algebraic is given in 9.4; the equivalence between
the splice and plumbing graphs is described in sections 20-22;
the splicing construction appears in section 8.

\subsection{The Alexander polynomial}\label{s0e}
 Assume that $K$ {\em is a homologically trivial oriented knot} in $M$. Let
$V$ be the Alexander matrix of $K\subset M$, and $V^*$ its transposed (cf.
\cite{Lescop}, page 26).

In the literature one can find different {\em normalizations} of the Alexander
``polynomial''. The most convenient for us, which makes our formulae
the simplest possible, is
$$\Delta^{\no}_M(K)(t):=\det (t^{1/2}\,V-t^{-1/2}\, V^*).$$
In the surgery formula  \ref{s3} we will need
{\em Lescop's normalization}  \cite{Lescop}, in this article
denoted by $\Delta^L_M(K)(t)$.
They are related by the identity (cf.  \cite[2.3.13]{Lescop}):
$$\Delta_M^{\no}(K)(t)=\Delta^L_M(K)(t)/|H|.$$
Then (see e.g. \cite[2.3.1]{Lescop}) one has:
\begin{equation*}
\Delta^{\no}_M(K)(1)=\Delta^L_M(K)(1)/|H|= 1.
\tag{6}
\end{equation*}
We also prefer  to think about
the Alexander polynomial as a characteristic polynomial.
For this notice that $V$ is invertible over $\Q$, hence one can
define the ``monodromy operator''  $\M:=V^{-1}V^*$. Then set:
\begin{equation*}
\Delta_M(K)(t):=\det (I-t\M)= \det(V^{-1})\cdot t^r\cdot \Delta_M^{\no}(K)(t).
\tag{7}
\end{equation*}
If $(M,K)$ can be  represented by a negative definite plumbing graph,
then by a theorem of Grauert \cite{Grauert}
$(M,K)$ is algebraic, hence by Milnor fibration
theorem, it is fibrable. In this case, $\M$  is exactly the monodromy operator
acting on the first homology of the (Milnor) fiber. Moreover,
$\Delta_M(K)(t)$ can be computed from the plumbing graph by
A'Campo's theorem \cite{AC} as follows (see also \cite{EN}).
Assume that $K=g_u$ for some  $u\in {\mathcal V}$, then
\begin{equation*}
\frac{\Delta_M(g_u)(t)}{t-1}= \prod_{v\in{\mathcal V}}\
(t^{w_v(u)}-1)^{\bar{\delta}_v-2}.
\tag{8}
\end{equation*}
Notice that (6) guarantees that $\Delta_M(K)(1)=\det(V^{-1})$, from
(8) we get $\Delta_M(K)(1)>0$, and from the Wang exact sequence of the
fibration $|\Delta_M(K)(1)|=|H|$. Therefore:
\begin{equation*}
\Delta_M(K)(1)=\det(V)^{-1}= |H|.
\tag{9}
\end{equation*}
More generally, if $(M,K)$ has a negative definite
plumbing representation, and $K=g_u$ for some
$u$, then for any character $\chi\in \hat{H}$ we define
$\Delta_{M,\chi}(g_u)(t)$ via the identity
\begin{equation*}
\frac{\Delta_{M,\chi}(g_u)(t)}{t-1}:=
\prod_{v\in{\mathcal V}}\
(t^{w_v(u)}\chi(g_v)-1)^{\bar{\delta}_v-2},
\tag{10}
\end{equation*}
and we write
\begin{equation*}
\Delta^H_M(g_u)(t):=\frac{1}{|H|}\cdot
\sum_{\chi\in \hat{H}}\ \Delta_{M,\chi}(g_u)(t).
\tag{11}\end{equation*}
%Fix a character $\chi$ and vertex $v$ with $\chi(g_v)=1$.
%Then the number of adjacent vertices $v'$  of  $v$ with $\chi(g_{v'})\not=1$
%cannot be one (cf. \cite{NN}, A7). Using this fact, and by analyzing
%the  support of a non-trivial character $\chi$  (i.e.  the
%subgraph spanned by those vertices $v$ with  $\chi(g_v)=1$), one gets
%(cf. \cite{NN}, A7):
%\begin{equation*}
%\lim_{t\to 1} \Delta_{M,\chi}(g_u)(t)=
%\left\{\begin{array}{ll}
%|H| &\mbox{if $\chi=1$},\\
%0 & \mbox{if $\chi\not=1$}.\end{array} \right.
%\tag{11}
%\end{equation*}
%In particular,
%\begin{equation*}
%\lim_{t\to 1}\ \Delta^H_{M}(g_u)(t)=1.
%\tag{12}
%\end{equation*}

In section 6 we will need the following analog of (9) in the case when $K=g_u
\subset M$ is not homologically trivial (but $(M,K)$ has a negative definite
plumbing representation):
\begin{equation*}
\lim_{t\to 1}\ (t-1) \prod_{v\in{\mathcal V}}\
(t^{w_v(u)}-1)^{\bar{\delta}_v-2}=|H|/o(u).
\tag{12}
\end{equation*}
This follows e.g. from \cite{NN} A10(b). (In fact, $|H|/o(u)$ has the geometric
meaning of $|Tors\, H_1(M\setminus K)|$, and (12) can also be deduced from the
Wang exact sequence of the monodromy, similarly as above.)

\subsection{The Seiberg-Witten invariant}\label{sw} If $M$ is a rational
homology 3-sphere, then
the set $Spin^c(M)$ of the $spin^c$ structures of $M$  is a $H$-torsor.
If $M$ is the link of a normal surface singularity (or, equivalently,
if $M$ has a plumbing representation with a negative definite
intersection matrix), then $Spin^c(M)$ has a distinguished element
$\sigma_{can}$,
called the {\em canonical } $spin^c$ structure  (cf. \cite{NN}).

To describe the
Seiberg-Witten invariants one has to consider  an additional geometric
data belonging to  the space of  parameters
\[
\p=\{ u=(g,\eta);\quad g=\mbox{Riemann metric},\;\;\eta=\mbox{closed
two-form} \}.
\]
Then for each $spin^c$ structure $\si$ on $M$  one defines the
$(\si,g,\eta)$-{\em Seiberg-Witten monopoles}.
For a generic parameter $u$, the Seiberg-Witten invariant  $\ssw_M(\si,u)$
is the signed monopole count.
This integer depends on the choice of the parameter $u$ and thus
it is not a topological  invariant.   To obtain an invariant of $M$,  one
needs  to  alter this monopole count. The  additional
contribution is the {\em Kreck-Stolz} invariant $KS_M(\si,u)$
(associated with the data $(\si,u)$),
cf. \cite{Lim} (or see \cite{KrSt} for the original ``spin version'').
Then, by \cite{Chen1,Lim,MW}, the rational number
\[
\frac{1}{8}KS_M(\si, u)+\ssw_M(\si, u)
\]
is independent   of $u$ and thus  it is a topological invariant of
the pair $(M,\si)$. We denote this {\em modified Seiberg-Witten
invariant}  by \ $\ssw_M^0(\si)$.
In general, it is very difficult to compute $\ssw_M^0(\si)$ using this
definition. Fortunately, it has another realization as well.
For any $spin^c$ structure $\si$ on $M$, we denote by
$$\et_{M,\si}=\sum_{h\in H}\et_{M,\si}(h)\, h\in {\bQ}[H]$$
 the sign refined
{\em Reidemeister-Turaev torsion} associated with $\si$
(for its detailed description, see \cite{Tu5}).
It is convenient to think of $\et_{M,\si}$ as a function
$H\ra {\bQ}$ given by $h\mapsto \et_{M,\si}(h)$.
The  {\em augmentation map} $\aug: {\bQ}[H]\ra {\bQ}$
is defined by $\sum a_h\, h\mapsto \sum a_h$. It is known that
$\aug(\et_{M, \si})=0$.

Denote the Casson-Walker invariant of $M$ by $\lambda_W(M)$ \cite{walker}.
It is related with Lescop's normalization $\lambda(M)$
 \cite[\S 4.7]{Lescop} by $\lambda_W(M)=2\lambda(M)/|H|$.
Then by  a result of the second author
\cite{Nico5}, one has:
\begin{equation*}
\ssw^0_M(\si)=\et_{M,\si}(1)-\lambda_W(M)/2.
\tag{13}
\label{3.7.2}
\end{equation*}
Below we will present a formula for
$\et_{M,\si}$ in terms of Fourier transform. Recall that
a function $f: H \ra {\bC}$ and
its Fourier transform $\hat{f}:\hat{H}\ra {\bC}$ satisfy:
 \[
\quad \hat{f}(\chi)=\sum_{h\in H} f(h)\bar{\chi}(h);\ \
 f(h)=\frac{1}{|H|}\sum_{\chi\in \hat{H}}\hat{f}(\chi)\chi(h).
 \]
Here $\hat{H}$ denotes the Pontryagin dual of $H$ as above.
 Notice  that $\hat{f}(1)=\aug(f)$,
in particular $\hat{\et}_{M,\si}(1)=\aug(\et_{M, \si})=0$.
Therefore,
\begin{equation*}
\et_{M,\sigma}(1)=\frac{1}{|H|}\sum_{\chi\in \hat{H}\setminus \{1\}}
\ \hat{\et}_{M,\sigma}(\chi).
\tag{14}
\end{equation*}
Now, assume that $M$ is represented by a
{\em negative definite plumbing graph}.
Fix a non-trivial character $\chi\in \hat{H}\setminus \{1\}$ and an arbitrary
vertex $u\in {\mathcal V}$ with $\chi(g_u)\not=1$. Then set
\begin{equation*}
\hat{P}_{M,\chi,u}(t):= \prod_{v\in{\mathcal V}}\
(t^{w_v(u)}\chi(g_v)-1)^{\delta_v-2},
\tag{15}\end{equation*}
where $t\in \C$ is a free variable. Then, by \cite{NN} (5.8),
the Fourier transform $\hat{\et}_{M,\si_{can}}$ of \ $\et_{M,\si_{can}}$
is given by
\begin{equation*}
\hat{\et}_{M,\si_{can}}(\bar{\chi})=
\lim_{t\to 1}\, \hat{P}_{M,\chi,u}(t).
\tag{16}
\end{equation*}
This  limit is independent of the choice of $u$, as long as
$\chi(g_u)\not=1$. In fact, by \cite{NN}, even if
$\chi(g_u)=1$, but $u$ is adjacent to some vertex $v$ with $\chi(g_v)\not=1$,
$u$ does the same job.

\section{Some general splicing formulae}

\subsection{The splicing data}\label{s1} We will consider the following
geometric  situation. We start with two oriented 3-manifolds $M_1$ and $M_2$,
both {\em rational homology spheres}. For $i=1,2$, we fix  an oriented  knot
$K_i$ in $M_i$, and we use the notations of \ref{s0} with the corresponding
indices $i=1,2$.
In this article we will consider a particular splicing, which is motivated
by the geometry of the suspension singularities. The more general case will
be treated in a forthcoming paper.

On the pair $(M_2,K_2)$ we impose no additional restrictions. But,
for $i=1$, we will consider the following {\em working assumption}:

\vspace{2mm}

\begin{description}
\item[WA1]
{\em Assume that $o_1=1$, i.e. $K_1 $ is
homologically trivial in $M_1$.  Moreover, we fix the parallel $\ell_1$
exactly as  the longitude $\lambda_1$. Evidently,  $k_1=0$. }
\end{description}

\vspace{2mm}

\noindent
Finally, by splicing,  we define a  3-manifold $M$  (for details, see e.g.
\cite{fujita}):
$$M=M_1\setminus T(K_1)^\circ \coprod_A M_2\setminus T(K_2)^\circ,$$
where $A$ is an identification of $\partial T(K_2)$ with $-\partial T(K_1)$
determined by
\begin{equation*}
A(m_2)=\lambda_1, \ \ \mbox{ and} \ \  A(\ell_2)=m_1.
\tag{1}\end{equation*}

\subsection{The closures $\overline{M}_i$}\label{s1b}
Once the  splicing data is fixed, one can consider the {\em closures}
$\overline{M}_i$ of $M_i\setminus T(K_i)^\circ $ ($i=1,2$) with respect to
$A$ (cf. \cite{BN} or \cite{fujita}) by the following Dehn fillings:
$$\overline{M}_2=(M_2\setminus  T(K_2)^\circ )(A^{-1}(y_1)), \ \
\overline{M}_1=(M_1\setminus  T(K_1)^\circ )(A(y_2)),$$
where $\delta_iy_i:=\lambda_i$ ($i=1,2$).
Using (1)  one has
$A^{-1}(y_1)=m_2$, hence
$$\overline{M}_2=M_2.$$
Moreover, $A(y_2)=A((o_2\ell_2+k_2m_2)/\delta_2)=(o_2m_1+k_2\lambda_1
)/\delta_2$, hence:
$$\overline{M}_1=(M_1\setminus T(K_1)^\circ)(\mu), \ \ \mbox{where}\ \
\mu:= (o_2m_1+k_2\lambda_1)/\delta_2.$$
In fact, $\overline{M}_1$ can be represented as a $(p,q)$-surgery of
$M_1$ along $K_1$. The integers $(p,q)$ can be determined as in \cite{Lescop},
page 8: $\mu$ is homologous to $qK_1$ in $T(K_1)$, hence
$q=k_2/\delta_2$. Moreover, $ p=Lk_{M_1}(\mu,K_1)$, which via \ref{s0c}(1)
equals
$\lag(o_2m_1+k_2\lambda_1)/\delta_2, \lambda_1\rag=o_2/\delta_2$.
Therefore,
\begin{equation*}
\overline{M}_1=M_1(K_1, p/q)=M_1(K_1,o_2/k_2).\tag{2}
\end{equation*}

\subsection{Fujita's splicing formula for the Casson-Walker
invariant.}\label{s2} Using the above expressions for the closures,
(1.1) of \cite{fujita},
in the case of the above splicing (with $A=f^{-1}$),  reads as
%$$\lambda_W(M)=\lambda_W(\overline{M}_1)+\lambda_W(\overline{M}_2)+
%\bms(-u_2,o_2/\delta_2).$$
%Using the identity \ref{s0}(1) one gets that $-u_2=(k_2/\delta_2)^{-1},
%\mmod ( o_2/\delta_2)$, hence $\bms(-u_2,o_2/\delta_2)=\bms(k_2/\delta_2,
%o_2/\delta_2)=\bms(k_2,o_2)$. Therefore,
\begin{equation*}
\lambda_W(M)=\lambda_W(M_2)+\lambda_W(M_1(K_1,o_2/k_2))+\bms(k_2,o_2).
\tag{3}\end{equation*}
Additionally, if  we assume that $K_2$ is  homologically trivial
in $M_2$ (i.e.  $o_2=1$), and we fix $\ell_2$ as $\lambda_2$ (i.e. $k_2=0$),
then (3) transforms into:
\begin{equation*}
\lambda_W(M)=\lambda_W(M_1)+\lambda_W(M_2).
\tag{4}\end{equation*}

\subsection{Walker-Lescop  surgery formula}\label{s3} Now, we will analyze
the manifold $M_1(K_1,p/q)$ obtained by $p/q$-surgery,
where $p=o_2/\delta_2>0$ and $q=k_2/\delta_2$
(not necessarily positive). First notice (cf.  \cite[1.3.4]{Lescop}) that
$|H_1(M_1(K_1,p/q))|=p\cdot |H_1(M_1)|$.
Using this, the surgery formula (T2) from \cite{Lescop}, page 13,
and the identification
$\lambda_W(\cdot)=2\lambda(\cdot)/|H_1(\cdot,\Z)|$,  one gets:
\begin{equation*}
\lambda_W(M_1(K_1,p/q))=\lambda_W(M_1)+Cor,
\tag{5}
\end{equation*}
where the correction term $Cor$ is
$$Cor:= \frac{q}{p}\cdot \frac{\Delta^L_{M_1}(K_1)''(1)}{|H_1(M_1)|}-
\frac{p^2+1+q^2}{12pq}+\sign(q)\Big( \frac{1}{4}+ \bms(p,q)\Big).$$
Using (3), (5) and the reciprocity law  of the Dedekind sums (for $p>0$):
$$ \bms(q,p)+\sign(q)\bms(p,q)=-\frac{\sign(q)}{4}+\frac{p^2+1+q^2}{12pq},$$
one gets the following formula:

\subsection{Theorem. The splicing formula for Casson-Walker
invariant.}\label{s4} {\em Consider a splicing manifold $M$
characterized by the
data  described in \ref{s1}. Then:
$$\lambda_W(M)=\lambda_W(M_1)+ \lambda_W(M_2)+\frac{k_2}{o_2}\cdot
\Delta^{\no}_{M_1}(K_1)''(1).$$}

\subsection{The splicing property of the group $H_1(M,\Z)$.}\label{s5} In the next
paragraphs we analyze  the behavior of $H_1(\cdot ,\Z)$ under the
splicing construction \ref{s1}.

First notice that by excision, for any $q$, one has:
$$H_q(M,M_2\setminus T(K_2)^\circ )=H_q(M_1\setminus  T(K_1)^\circ ,
\partial T(K_1))=H_q(M_1,K_1).$$
Therefore, the long exact
sequence of the pair $(M,M_2\setminus T(K_2)^\circ)$ reads as:
$$0\to H_2(M_1,K_1)\stackrel{\partial_1}{\longrightarrow}
H_1(M_2\setminus T(K_2)^\circ )\to H_1(M)\to H_1(M_1,K_1)\to 0.$$
Using the isomorphisms \ref{s0c}(3), $\partial_1$ can be identified with $
\partial_1(1_\Z)=m_2$, hence coker$\partial_1=H_1(M_2)$, cf. \ref{s0c}(2).
Therefore,  \ref{s0c}(3) gives the exact sequence
\begin{equation*}
0\to H_1(M_2)\stackrel{i}{\longrightarrow}
H_1(M)\stackrel{p}{\longrightarrow} H_1(M_1)\to 0.
\tag{6}
\end{equation*}
This exact sequence splits. Indeed, let $\bar{s}$ be the composition of
$s_1:H_1(M_1)\to H_1(M_1\setminus T(K_1))$ (cf. \ref{s0c})
and $H_1(M_1\setminus T(K_1))\to
H_1(M)$ (induced by the inclusion). Then $p\circ \bar{s}=1$. In particular,
\begin{equation*}
H_1(M)=im(i)\oplus im(\bar{s})\approx H_1(M_2)\times H_1(M_1).
\tag{7}
\end{equation*}
Notice that any $[K]\in H_1(M_1)$ can be represented (via $s_1$) by
a representative $K$ in $M_1\setminus T(K_1)$ providing a class in
$Tors\, H_1(M_1\setminus T(K_1))$. Write
$o$ for its order, and take a Seifert surface $F$, sitting in
$M_1\setminus T(K_1)$  with $\partial F=oK$.
If $L\subset M_1\setminus T(K_1)$ with $L\cap K=\varnothing$
 then obviously $Lk_M(K,L)=Lk_{M_1}(K,L)$.
Moreover, since $F$ has no intersection
points with any curve  $L\in M_2\setminus T(K_2)$, for such an $L$ one gets:
\begin{equation*}
Lk_M(K,L)=0, \ \ \mbox{hence} \ \ b_M(im(\bar{s}), im (i))=0.
\tag{8}
\end{equation*}
By a similar argument,
\begin{equation*}
b_M(\bar{s}(x),y)=b_{M_1}(x,p(y))\ \ \mbox{
for any $x\in H_1(M_1)$ and $y\in H_1(M)$.}
\tag{9}\end{equation*}
The point is that in the next sections we need $Lk_M(K,\cdot)$
for general $K\subset M_1\setminus K_1$ which
is not a torsion element in $H_1(M_1\setminus K_1)$.

To compute this linking number  consider  another oriented knot $L\subset M_1\setminus K_1$ with $K\cap L=\varnothing$.
Our goal is to compare $Lk_{M_1}(K,L)$ and  $Lk_M(K,L)$.
Assume that the order of the class $K$ in $H_1(M_1)$ is $o$.
Let $F_{oK}$  be a Seifert surface in $M_1$ with $\partial F_{oK}=oK$ which
intersects $L$ and $K_1$ transversally. It is clear that $F_{oK}$ intersects
$L$  exactly in $o\cdot Lk_{M_1}(K,L)$  points (counted with sign).
On the other hand, it intersects
$K_1$ in $o\cdot Lk_{M_1}(K,K_1)$ points. For each
intersection point with sign $\eps=\pm1 $, we cut out from $F_{oK}$ the disc $F_{oK}\cap T(K_1)$, whose orientation depends on $\eps$. Its boundary is $\eps m_1$ which, by the splicing
identification, corresponds to $\eps\ell_2$ in $M_2\setminus T(K_2)$.
Using rational coefficients, $\ell_2=(1/o_2)\lambda_2-(k_2/o_2)m_2$.
Some multiple of $\lambda_2$  can be extended to a surface in $M_2\setminus
K_2$ which clearly has no intersection with $L$. $m_2$ by splicing is
identified with $\lambda_1$ which has a Seifert surface in $M_1\setminus K_1$
which intersect $L$ in $Lk_{M_1}(K_1,L)$ points. This shows that
for $K,L\subset M_1\setminus T(K_1)$:
\begin{equation*}
Lk_M(K,L)=Lk_{M_1}(K,L)-Lk_{M_1}(K,K_1)\cdot Lk_{M_1}(L,K_1)\cdot k_2/o_2.
\tag{10}
\end{equation*}
Assume now that $K\in M_2\setminus T(K_2)$, and let $o$ be the order of its
homology class in $H_1(M_2)$. Assume that the
Seifert surface $F$ of $oK$ intersects $K_2$ transversally in $t_2$ points.
Then $F\cap \partial T(K_2)=t_2m_2$, and after the splicing identification
this becomes $t_2\lambda_1$. Let $F_1$ be the Seifert surface of
$\lambda_1$ in $M_1\setminus  T(K_1)^\circ$. Then (after some natural
identifications) $F\setminus T(K_2)^\circ \coprod t_2F_1$ is a Seifert
surface of $oK$ in $M$. Therefore, for any $L\subset  M_1\setminus T(K_1)$
one has:
\begin{equation*}
Lk_M(K,L)=Lk_{M_2}(K,K_2)\cdot Lk_{M_1}(K_1,L)\ \ \mbox{for
$K\subset M_2\setminus T(K_2)$ and $L\subset M_1\setminus T(K_1)$}.
\tag{11}
\end{equation*}
By a  similar  argument:
\begin{equation*}
Lk_M(K,L)=Lk_{M_2}(K,L)\ \ \mbox{for
$K, L\subset M_2\setminus T(K_2)$}.
\tag{12}
\end{equation*}

\subsection{Splicing plumbing manifolds}\label{s6} Some of the
results of this article about Reidemeister-Turaev torsion can be formulated
and proved in the context of general (rational homology sphere)
3-manifolds, and arbitrary $spin^c$ structures. Nevertheless, in this
article we are mainly interested in algebraic links, therefore
we restrict ourselves to plumbed manifolds.

This can be formulated in the second {\em working assumption}:

\vspace{2mm}

\begin{description}
\item[WA2] {\em
$(M_i,K_i) \ (i=1,2)$ can be represented by a
negative definite plumbing graph. Moreover, if $M$ is the result of
the splicing (satisfying  WA1, cf.  \ref{s1}),
then $M$ can also be represented by a negative
definite plumbing graph.}
\end{description}

\vspace{2mm}

\noindent  Assume that the plumbing graphs
$\Gamma(M_1,K_1)$ and $\Gamma(M_2,K_2)$  have the following schematic
form (with $g_{v_1}=K_1$ and $g_{v_2}=K_2$):

\begin{picture}(400,40)(-100,30)
\put(20,40){\framebox(40,20)}
\put(70,50){\circle*{4}}
\put(70,50){\line(-2,1){15}}
\put(70,50){\line(-2,-1){15}}
\put(55,52){\makebox(0,0){$\vdots$}}
\put(70,50){\vector(1,0){20}}
\put(70,60){\makebox(0,0){$v_1$}}
\put(160,40){\framebox(40,20)}
\put(150,50){\circle*{4}}
\put(150,50){\line(2,1){15}}
\put(150,50){\line(2,-1){15}}
\put(165,52){\makebox(0,0){$\vdots$}}
\put(150,50){\vector(-1,0){20}}
\put(150,60){\makebox(0,0){$v_2$}}
\end{picture}

\noindent
Then it is not difficult to see that (a possible) plumbing graph $\Gamma(M)$
for $M$ has the following form (where $v_1$ and $v_2$ are connected by a
string):

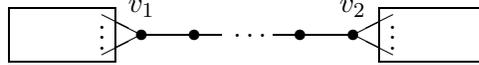
\begin{picture}(400,40)(-100,30)
\put(20,40){\framebox(40,20)}
\put(70,50){\circle*{4}}
\put(90,50){\circle*{4}}
\put(70,50){\line(-2,1){15}}
\put(70,50){\line(-2,-1){15}}
\put(55,52){\makebox(0,0){$\vdots$}}
\put(112,50){\makebox(0,0){$\dots$}}
\put(70,50){\line(1,0){30}}
\put(70,60){\makebox(0,0){$v_1$}}
\put(160,40){\framebox(40,20)}
\put(150,50){\circle*{4}}
\put(130,50){\circle*{4}}
\put(150,50){\line(2,1){15}}
\put(150,50){\line(2,-1){15}}
\put(165,52){\makebox(0,0){$\vdots$}}
\put(150,50){\line(-1,0){30}}
\put(150,60){\makebox(0,0){$v_2$}}
\end{picture}

If ${\mathcal V}_i$ ($i=1,2$), respectively ${\mathcal V}$,
 represent the set of vertices of $\Gamma(M_i)$, resp. of
$\Gamma(M)$, then ${\mathcal V}={\mathcal V}_1\cup {\mathcal V}_2$
modulo some vertices with $\delta=2$. (In particular, in any formula like
\ref{sw}(15), ${\mathcal V}$ behaves as the union
${\mathcal V}_1\cup {\mathcal V}_2$.)

\subsection{Remark.}\label{s6r} If one wants to compute $\et_{M,\sigma_{can}}$
for such plumbed manifolds, then one can apply \ref{sw}(15) and (16).
For this, one has to analyze the supports of the characters, and the
corresponding weights $w_v(u)$. These weights are closely related
with the corresponding linking numbers $Lk_M(g_u,g_v)$ (cf. \ref{s0d} (5)),
hence the relations \ref{s5}(10)-(11)-(12) are crucial.
For characters of type $\chi\in \hat{p}(H_1(M_1)\hat{}\,)$
(cf. \ref{s5}(6) or the next proof)
$u\in {\mathcal V}_1$, hence \ref{s5}(10) should be applied for any
$v\in {\mathcal V}_1$. But this is
rather  unpleasant due to the term $Lk_{M_1}(g_u,K_1)\cdot
Lk_{M_1}(g_v,K_1)\cdot k_2/o_2$. The description is more transparent if
either $H_1(M_1)=0$ or $k_2=0$.

Therefore, we will consider first these particular cases only.
They, as guiding examples,
 already contain all the illuminating information and principles
we need to proceed.
For the link of $\{f(x,y)+z^n=0\}$,
the splicing formula will be made very explicit  in \ref{splet}
(based an a detailed and complete classification of the characters
and the regularization terms $\hat{P}$, which is rather involved).

\subsection{Theorem. Some splicing formulae for the Reidemeister-Turaev
torsion}\label{s7} {\em Assume that $M$ satisfy WA1 (\ref{s1}) and WA2
(\ref{s6}). Then the following hold:

\noindent  (A) Assume that $K_2\subset M_2$ is also homologically
trivial (i.e. $o_2=1$), and $\ell_2=\lambda_2$ (i.e. $k_2=0$).
Then
$$\et_{M,\sigma_{can}}(1)=\et_{M_1,\sigma_{can}}(1)+
\et_{M_2,\sigma_{can}}(1).$$

\noindent (B) Assume that $M_1$ is an integral homology sphere (i.e.
$H_1(M_1)=0$). Then
$$\et_{M,\sigma_{can}}(1)= \sum_{\chi_2\in H_1(M_2)\hat{}\,\setminus\{1\}}\
\frac{\hat{\et}_{M_2,\sigma_{can}}(\bar{\chi}_2)}{|H_1(M_2)|}\cdot
\Delta_{M_1}(K_1)(\chi_2(K_2)).$$
In particular, if $K_2\subset M_2$ is homologically trivial, then $\chi_2(K_2)
=1$ for any $\chi_2$, hence by \ref{s0e}(9)  one gets
$$\et_{M,\sigma_{can}}(1)=\et_{M_2,\sigma_{can}}(1)
\ \ \ (\mbox{and evidently \ $\et_{M_1,\sigma_{can}}(1)=0$}).$$
This is true for any choice of $\ell_2$, i.e. even if $k_2$ is non-zero.}

\begin{proof}
The theorem is a consequence of \ref{sw}(14)-(15)-(16)
and \ref{s5}. For this,  we have to analyze the characters
$\chi$ of $H_1(M)$.  The dual of the exact sequence \ref{s5}(6) is
$$0\to H_1(M_1)\hat{}
\stackrel{\hat{p}}{\longrightarrow}
H_1(M)\hat{} \stackrel{\hat{i}}{\longrightarrow }
H_1(M_2)\hat{}\to 0.$$
First, consider a character $\chi$ of $H_1(M)$ of the form $\chi=\hat{p}
(\chi_1)$ for some $\chi_1\in H_1(M_1)\hat{}$. Since any $\chi_1\in H_1
(M_1)\hat{}$ can be represented as $b_{M_1}(x,\cdot)$ for some $x\in H_1(M_1)$
(cf. \ref{s0c}), and $\hat{p}(b_{M_1}(x,\cdot ))=b_M(\bar{s}(x),\cdot)$
(cf. \ref{s5}(9)), property \ref{s5}(8) guarantees that $\chi(g_v)=1$
for any $v\in {\mathcal V}_2$.
In particular, for $\chi=\hat{p}(\chi_1)$ with $\chi_1\in H_1(M_1)\hat{}
\setminus \{1\}$,  and for some $u\in {\mathcal V_1}$ with $\chi_1(g_u)\not=1$
(which works for $M$ as well), one gets:
$$\hat{P}_{M,\chi,u}(t)=\prod_{v\in {\mathcal V}_1}\
(t^{w_v(u)}\chi_1(g_v)-1)^{\bar{\delta}_v-2}\cdot
\prod_{v\in {\mathcal V}_2}\
(t^{w_v(u)}-1)^{\bar{\delta}_v-2}.$$
Here, for $v\in {\mathcal V}_i$, $\bar{\delta}_v$ means the number of adjacent
edges of $v$ in $\Gamma(M_i,K_i)$ ($i=1,2$), cf. \ref{s0d}.

By \ref{s5}(11), for any $v\in {\mathcal V}_2$ one has
$$w_v(u)=o(u)Lk_M(g_u,g_v)=o(u)Lk_{M_1}(g_u,g_{v_1})\cdot Lk_{M_2}(g_{v_2},
g_v),$$
hence by \ref{s0e}(8)
$$\hat{P}_{M,\chi,u}(t)=\hat{P}_{M_1,\chi_1,u}(t) \cdot
\Delta_{M_2}(g_{v_2})(t^{o(u)Lk_{M_1}(g_u,g_{v_1})}).$$
Taking the limit $t\to 1$, and using \ref{s0e}(9), one gets:
\begin{equation*}
\hat{\et}_{M,\sigma_{can}}(\bar{\chi})=
\hat{\et}_{M_1,\sigma_{can}}(\bar{\chi}_1)\cdot |H_1(M_2)|.
\tag{13}
\end{equation*}
Now, we prove (A). In this case $M_1$ and $M_2$ are symmetric, hence
there is a similar term as in (13) for characters $\chi_2\in H_1(M_2)\hat{}$.

On the other hand, if $\chi=\chi_1\chi_2$ for two non-trivial characters
$\chi_i\in H_1(M_i)\hat{}$ ($i=1,2$), then one can show  that
$\hat{P}_{M,\chi,u}(t)$ has a root (of multiplicity at least
two) at $t=1$, hence $\hat{\et}_{M,\sigma_{can}}(\bar{\chi})=0$. Now, use
\ref{sw}(14) and
$|H_1(M)|=|H_1(M_1)|\cdot |H_1(M_2)|$ (cf. \ref{s5}(7)).

For (B), fix a non-trivial character $\chi_2\in H_1(M_2)\hat{}$.
The relations in \ref{s5} guarantee that if we take
$u\in {\mathcal V}_2$ with $\chi_2(g_u)\not=1$ (considered as a property
of $M_2$) then $\chi(g_u)\not=1$ as well for $\chi= \hat{i}^{-1}(\chi_2)$.
Moreover, for any $v\in {\mathcal V}_1$ one has:
$\chi(g_v)=\chi_2(K_2)^{Lk_{M_1}(K_1,g_v)}$ (cf. \ref{s5}(11)). Therefore:
$$\hat{P}_{M,\chi,u}(t)=\prod_{v\in {\mathcal V}_1}\
(t^{w_v(u)}\chi(g_v)-1)^{\bar{\delta}_v-2}\cdot
\prod_{v\in {\mathcal V}_2}\
(t^{w_v(u)}\chi(g_v)-1)^{\bar{\delta}_v-2}$$
%$$=\hat{P}_{M_1,\chi_1,u}(t) \cdot
%(t^{o(u)Lk_{M_1}(g_u,g_{v_1})}\chi_1(g_{v_1})-1)\cdot
%\prod_{v\in {\mathcal V}_2} (t^{o(u)Lk_{M_1}(g_u,g_{v_1})Lk_{M_2}(g_v,g_{v_2})}
%\chi_1(K_1)^{Lk_{M_2}(g_v,g_{v_2})}\chi'(g_v)-1)^{\bar{\delta}-2}$$
$$=\hat{P}_{M_2,\chi_2,u}(t) \cdot
\Delta_{M_1}(g_{v_1})(t^{ o(u)Lk_{M_2}(g_u,g_{v_2})}\chi_2(K_2)).$$
\end{proof}

\subsection{Remarks.}\label{s9} \
(1) A similar proof provides the following
formula as well (which  will be not used later).
Assume that $M$ satisfies WA1 and WA2, and $k_2=0$. Then
$$\et_{M,\sigma_{can}}(1)=\et_{M_1,\sigma_{can}}(1)+
 \sum_{\chi_2\in H_1(M_2)\hat{}\,\setminus\{1\}}\
\frac{\hat{\et}_{M_2,\sigma_{can}}(\bar{\chi}_2)}{|H_1(M_2)|}\cdot
\Delta^H_{M_1}(K_1)(\chi_2(K_2)).$$
%(If $k_2\not=0$ and $H_1(M_1)\not=0$ then the first summand in the right hand
%side  should be substantially modified.)

(2) The obstruction term (see also \ref{def1}(1)) which measures the
non-additivity of the Casson-Walker invariant (under the splicing assumption
WA1) is given by
$\frac{k_2}{o_2}\cdot  \Delta^{\no}_{M_1}(K_1)''(1)$ (cf. \ref{s4}).
On the other hand, if $H_1(M_1)=0$, then the obstruction term for
the non-additivity of the
Reidemeister-Turaev torsion (associated with $\sigma_{can}$) is
$$\sum_{\chi_2\in H_1(M_2)\hat{}\, \setminus\{1\}}\
\frac{\hat{\et}_{M_2,\sigma_{can}}(\bar{\chi}_2)}{|H_1(M_2)|}\cdot
\Big(\ \Delta_{M_1}(K_1)(\chi_2(K_2)) -1\Big).$$
Notice that they ``look rather different'' (even under  this extra-assumption
$H_1(M_1)=0$). In fact, even their  nature are different:
the first depends essentially on the choice of the parallel $\ell_2$
(see the coefficient $k_2$ in its expression), while the second not.
 In particular,  one cannot really hope (in general)
for the additivity of the modified Seiberg-Witten invariant.

Therefore, it is really remarkable and surprising, that in some of the
geometric situations discussed in the next sections, the
modified Seiberg-Witten invariant is additive
(though the invariants
$\lambda_W$ and $\et_{\sigma_{can}}(1)$ are not  additive,
their obstruction terms cancel each  other).

\section{The basic topological example}

\subsection{}\label{g1} Recall that in our main applications (for algebraic
singularities)
the involved 3-manifolds are plumbed manifolds. In particular, they can be
constructed inductively from Seifert manifolds by splicing (cf. also with
\ref{spl}). The present section has a double role.  First, we
work out  explicitly the splicing results obtained in the previous section
for the  case when $M_2$ is a Seifert manifold (with
special Seifert invariants). On the other hand, the detailed study
of this splicing formulae provides us a better understanding of the
subtlety of the behavior of the (modified) Seiberg-Witten  invariant
with respect to splicing and cyclic covers.
They will be formulated in some
  ``almost--additivity'' properties, where the non-additivity will be
characterized by a new invariant ${\mathcal D}$ constructed
from the Alexander polynomial of $(M_1,K_1)$.

\subsection{The splicing component $M_2$}\label{g2} Assume that $M_2$
is the link
$$\Sigma=
\Sigma(p,a,n):=\{(x,y,z)\in \C^3:\ x^p+y^a+z^n=0, \ |x|^2+|y|^2+|z|^2=1\}$$
of the Brieskorn hypersurface singularity
$X_2:=\{(x,y,z)\in \C^3:\ x^p+y^a+z^n=0\}$,
where gcd$(n,a)=1$ and gcd$(p,a)=1$. Set $d:=$gcd$(n,p)$. Then
$M_2$ is a rational homology sphere with the following data
(for more details, and for a more complete list of the relevant
invariants, see e.g. \cite{NN}, section 6; cf also with  \cite{JN,NeR}).

\vspace{2mm}

\noindent (a) \ the Seifert invariants are: $n/d, \ p/d,\
a, \ a, \cdots, a$ ($a$ appearing  $d$ times, hence all together there are
$d+2$ special fibers); these numbers also give (up to a sign) the determinants
of the corresponding arms of the plumbing graph of $\Sigma$;

\vspace{2mm}

\noindent (b) \  the orbifold Euler characteristic is
$e=-d^2/(npa)$;

\vspace{2mm}

\noindent (c) \ $H_1(M_2)=\Z_a^{d-1}$ (cf. also with \ref{g3}(1) below);

\vspace{2mm}

\noindent (e) \ $\et_{M_2,\sigma_{can}}(1)=\frac{np}{24da}(d-1)(a^2-1)$
(see also ($*$) in the proof of \ref{g8}).

\subsection{The link $K_2$}\label{g3} In this section we will modify
slightly the construction of \ref{s1}: we will fix a link with $d$ connected
components (instead of a knot), and we will perform splicing along each
connected component.

This link $K_2\subset M_2$ is given by the equation $\{y=0\}$ in $M_2$.
In other words, $K_2$ is the union of the $d$ special (Seifert) fibers
corresponding to the  Seifert invariant $a$. The components of $K_2$
will be denoted by $K_2^{(i)}$ ($1=1,\ldots,d$), and their tubular
neighborhood by $T(K_2^{(i)})$ as in \ref{s0}.

The link-components $\{K_2^{(i)}\}_{i}$ generate $H_1(M_2)$
(see, e.g. \cite{NN}, section 6); in fact $H_1(M_2)$ has the following
presentation (written additively):
\begin{equation*}
H_1(M_2)=\lan [K_2^{(1)}],\ldots,
[K_2^{(d)}]:\
a[K_2^{(i)}]\ \mbox{for each $i$, and} \
[K_2^{(1)}]+\cdots +[K_2^{(d)}]=0\ran.
\tag{1}
\end{equation*}
Besides  $K_2$, there are two more {\em special} orbits in $M_2=\Sigma$, namely
 $Z:=\{z=0\}$ and $X:=\{x=0\}$. Moreover, let
$O$  be the generic fiber of the Seifert fibration of $\Sigma$
(i.e. $g_v$ associated with the central vertex $v$).
Then one has the following linking numbers
(use \ref{s0d}(5) and \cite{NN} (5.5)(1)):
\begin{equation*}
\begin{array}{ll}
(a)\ Lk_{M_2}(K_2^{(i)},K_2^{(j)})=np/(d^2a) & \mbox{ for any $i\not=j$};\\
(b)\ Lk_{M_2}(K_2^{(i)},O)=np/d^2; \ \ Lk_{M_2}(K_2^{(i)},Z)=p/d &
\mbox{ for any $i$};\\
(c)\ Lk_{M_2}(O,O)=npa/d^2 \ \ \mbox{and}\ \ Lk_{M_2}(O,Z)=pa/d;& \ \\
(d)\ Lk_{M_2}(X,Z)=a. & \
\end{array}
\tag{L\#}
\end{equation*}
Notice that $K_2\subset M_2$ is fibrable. Indeed, $K_2=\{y=0\}$
is the link associated with the algebraic germ $y:(X_2,0)\to (\C,0)$,
hence one can take its Milnor fibration. Let $F$ be the fiber with
$\partial F=K_2$ (equivalently, take any minimal Seifert surface $F$ with
$\partial F=K_2$). Then for each  $i=1,\ldots, d$, we define the {\em parallel}
$\ell_2^{(i)}$ in $\partial T(K_2^{(i)})$ by $F\cap \partial T(K_2^{(i)})$.

Let $\lambda_2^{(i)}$ be the {\em longitude} of $K_2^{(i)}\subset M_2$, and
consider the invariants $o_2^{(i)},\  k_2^{(i)}$,
%\ \delta_2^{(i)},\  x_2^{(i)},\ y_2^{(i)}$,
etc. as in \ref{s0} with the corresponding  sub- and superscripts added.

\subsection{Lemma.}\label{g4} {\em For each $i=1,\ldots, d$ one has:
$$o_2^{(i)}=a \ \ \mbox{and}\ \
k_2^{(i)}=\frac{np(d-1)}{d^2}.$$}
%In particular, $\delta_1^{(i)}=$gcd$(a,d-1)$.}
\begin{proof}
The first identity is clear (cf. \ref{g3}(1)).
For the second, notice that (cf.
\ref{s0c}(1))
\[
-k_2^{(i)}/o_2^{(i)}=\lan \ell_2^{(i)},
\lambda_2^{(i)}\ran/o_2^{(i)}=Lk_{M_2}(\ell_2^{(i)},K_2^{(i)}).
\]
  Moreover,
\[
Lk_{M_2}(\sum_j\ell_2^{(j)},K_2^{(i)})=0,
\]
 hence
 \[
 k_2^{(i)}/a=\sum_{j\not=i}Lk_{M_2}(K_2^{(j)},K_2^{(i)}).
 \]
  Then use \ref{g3}(L\#a).
%
%let us represent ${\mathcal O}$
%in $\partial T(K_1^{(i)})$ as a linear combination of $\ell_1^{(i)}$ and
%$m_1^{(i)}$. The coefficient of $\ell_1^{(i)}$ is clearly $a$ (from the
%definition of the Seifert invariants). The coefficient of $m_1^{(i)}$
%is $\lan {\mathcal O}, \ell_1^{(i)}\ran =Lk_{M_1}({\mathcal O}, K_1)=
%dnp/d^2$ (by \ref{g2}(f)). Therefore,
%$${\mathcal O}=a \ell_1^{(i)} +\frac{np}{d} m_1^{(i)}\ \ \mbox{in
%$\partial T(K_1^{(i)})$}.$$
%Since
%$\lambda_1^{(i)}=a \ell_1^{(i)} + k_1^{(i)} m_1^{(i)}$ (cf. \ref{s0}),
% using \ref{s0c}(2) one gets $Lk_{M_1}({\mathcal O}, K_1^{(i)})=
%\lan {\mathcal O}, \lambda_1^{(i)}\ran/o_1^{(i)}=np/d-k_1^{(i)}$.
%But, by \ref{g2}(f), this number equals $np/d^2$.
\end{proof}

\subsection{The manifold $M$.}\label{g5} Next, we consider $d$ manifolds
$M_1^{(i)}$ with knots $K_1^{(i)}\subset M_1^{(i)}$ ($i=1,\ldots, d$),
each satisfying the assumption  WA1 (i.e. $o_1^{(i)}=1$,
$\ell_1^{(i)}=\lambda_1^{(i)}$, and $k_1^{(i)}=0$, cf. \ref{s1}).
Then, for each $i=1,\ldots , d$, we consider the splicing identification
of $\partial T(K_2^{(i)})$ with  $-\partial T(K_1^{(i)})$ (similarly as
in \ref{s1}(1)):
$$A^{(i)}(m_2^{(i)})=\lambda_1^{(i)} \ \ \mbox{and} \ \ \
A^{(i)}(\ell_2^{(i)})=m_1^{(i)} .$$
Schematically:

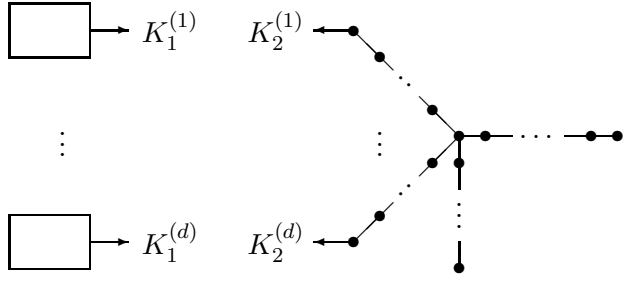
\begin{picture}(-400,110)(-400,0)
\put(-100,60){\circle*{4}}
\put(-90,60){\circle*{4}}
\put(-50,60){\circle*{4}}
\put(-40,60){\circle*{4}}
\put(-100,60){\line(0,-1){20}}
\put(-100,60){\line(-1,1){15}}
\put(-100,60){\line(-1,-1){15}}
\put(-100,10){\line(0,1){10}}
\put(-100,60){\line(1,0){20}}
\put(-40,60){\line(-1,0){20}}
\put(-70,60){\makebox(0,0){$\ldots$}}
\put(-100,33){\makebox(0,0){$\vdots$}}
\put(-130,60){\makebox(0,0){$\vdots$}}
\put(-250,60){\makebox(0,0){$\vdots$}}
\put(-100,50){\circle*{4}}
\put(-100,10){\circle*{4}}
%\put(-100,0){\circle*{4}}
\put(-119,41){\makebox(0,0){$\cdot$}}
\put(-122,38){\makebox(0,0){$\cdot$}}
\put(-119,79){\makebox(0,0){$\cdot$}}
\put(-122,82){\makebox(0,0){$\cdot$}}
\put(-110,70){\circle*{4}}
\put(-110,50){\circle*{4}}
\put(-130,90){\circle*{4}}
\put(-130,30){\circle*{4}}
\put(-140,100){\circle*{4}}
\put(-140,20){\circle*{4}}
\put(-140,100){\line(1,-1){15}}
\put(-140,20){\line(1,1){15}}
\put(-140,100){\vector(-1,0){15}}
\put(-140,20){\vector(-1,0){15}}
\put(-180,100){\makebox(0,0)[l]{$K_2^{(1)}$}}
\put(-180,20){\makebox(0,0)[l]{$K_2^{(d)}$}}
\put(-240,100){\vector(1,0){15}}
\put(-240,20){\vector(1,0){15}}
\put(-220,100){\makebox(0,0)[l]{$K_1^{(1)}$}}
\put(-220,20){\makebox(0,0)[l]{$K_1^{(d)}$}}
\put(-270,10){\framebox(30,20){}}
\put(-270,90){\framebox(30,20){}}
\end{picture}

\vspace{2mm}

\noindent {\em
In the sequel we denote by  WA1' the  assumption which guarantees
that  the manifold $M$ is constructed by this splicing procedure.
Moreover,  WA2' guarantees that all the involved 3-manifolds
have negative definite plumbing representations. }

\vspace{2mm}

\noindent Here some comments are in order.

\vspace{2mm}

\noindent (1) Assume that for some $i$, $M_1^{(i)}=S^3$, and $K_1^{(i)}$
is the unknot $S^1$ in $ S^3$. Then performing splicing along $K_2^{(i)}$
with $(S^3,S^1)$ is equivalent to put back $T(K_2^{(i)})$ unmodified,
hence it has no effect. In this case, one also has
$ \Delta_{S^3}(S^1)(t)\equiv \Delta^{\no}_{S^3}(S^1)(t)\equiv  1$.  \\
(2) Assume that we already had performed the splicing along the link-components
$K_2^{(i)}$ for $i\leq k-1$, but not along the other ones. Let us denote
the result of this partial modification by $M^{(k-1)}$. Consider $K_2^{(k)}$
in $M^{(k-1)}$ (in a natural way). Then all the invariants (e.g.
$o_2^{(k)}$,  $\lambda_2^{(k)}$,  $k_2^{(k)}$, etc.) associated with
 $K_2^{(k)}$  in $M_2$ or in $M^{(k-1)}$ are the same. (This follows from the
discussion in \ref{s5}, and basically, it is a consequence of WA1'.)

In particular, performing splicing at place  $i$ does not effect the splicing
data of the place $j$ ($j\not=i$). Therefore, using induction,
the computation of the invariants can be reduced easily
to the formulae established in the previous  section.

\subsection{Definitions/Notations.}\label{def1} \

\noindent
(1) In order to simplify the exposition, for any 3-manifold invariant
${\mathcal I}$,   we write
$${\mathcal O}({\mathcal I}):={\mathcal I}(M)-
{\mathcal I}(M_2)-\sum_{i=1}^d{\mathcal I} (M_1^{(i)}),$$
for  the ``additivity obstruction''  of ${\mathcal I}$ (with respect to the
splicing construction WA1'). E.g., using \ref{s5} and \ref{g5}(2)
one has  ${\mathcal O}(\log|H_1(\cdot)|)=0$.
Moreover, in all our Alexander invariant notations (e.g. in
$\Delta_{M_1^{(i)}}(K_1^{(i)})(t)$), we omit the link $K_1^{(i)}$
(e.g. we simply write $\Delta_{M_1^{(i)}}(t)$).

When we will compare ${\mathcal O}(\et_{\cdot,\sigma_{can}})$ with
${\mathcal O}(\lambda_W(\cdot ))$,
the next terminology will be helpful.

\noindent  (2) For any set of integers $c_1, \ldots, c_r$ define
$${\mathcal D}(c_1,\ldots, c_r):=
\sum_{i,j=1}^r \, c_ic_j\min(i,j)-\sum_{i=1}^r\, i\, c_i.$$

\noindent
(3) Define ${\mathcal D}(\Delta^{\no}(t))$ by ${\mathcal D}(c_1,\ldots, c_r)$
for any symmetric polynomial
\begin{equation*}
\Delta^{\no}(t)=1+\sum_{i=1}^r c_i(t^i+t^{-i}-2) \ \ (\mbox{for some
$c_i\in \Z$}).
%\tag{SymCoef}
\end{equation*}

\noindent (4)   A  set $\{c_i\}_{i\in I}$ ($I\subset \N$)
is called {\em alternating} if $c_i\in \{-1,0,+1\}$ for any $i\in I$; and
if $c_i\not=0$ then $c_i=(-1)^{n_i}$, where $n_i=\#\{j: j>i, \ \mbox{and}
\ c_j\not=0\}$.

\subsection{Corollary.}\label{g7}
{\em Assume that $M$ satisfies WA1'. Then }
$${\mathcal O}(\lambda_W)= \frac{np(d-1)}{ad^2}\cdot \sum_{i=1}^d
\, (\Delta^{\no}_{M_1^{(i)}})''(1).$$
\begin{proof} Use \ref{s4}, \ref{g4} and \ref{g5}(2). \end{proof}

\subsection{Corollary.}\label{g8}
{\em Assume that $M$ satisfies WA1' and WA2', and $M_1^{(i)}$ is an
integral homology sphere for any $i$.
Identify $\Z_a:=\{\xi\in \C: \xi^a=1\}$ and write $\Z_a^*:=\Z_a\setminus
\{1\}$.  Then }
$${\mathcal O}(\et_{\cdot ,\sigma_{can}}(1))=
\frac{np}{ad^2}\sum_{i\not=j}\sum_{\xi\in \Z_a^*}\
\frac{\Delta_{M_1^{(i)}}(\xi)\cdot \Delta_{M_1^{(j)}}(\bar{\xi})-1}
{(\xi-1)(\bar{\xi}-1)}.$$
\begin{proof} We recall first how one computes  the torsion for
the manifold $M_2$ (for a detailed  presentation, see \cite{NN}).
The point is (cf. also with the last sentence of \ref{sw}),
that for any character $\chi\in H_1(M_2)\hat{}\setminus \{1\}$,
one can choose the  central vertex of the star-shaped graph
for the vertex $u$ in order to  generate the weights in $\hat{P}$.
Then, by \cite{NN} or \ref{g3}(L\#), one gets:
$$\hat{P}_{M_2,\chi,u}(t)=\frac{(t^\alpha-1)^d}{(t^{d\alpha/n}-1)
(t^{d\alpha/p}-1)\prod_i\big(t^{\alpha/a}\chi(K_2^{(i)})-1\big)},$$
where $\alpha:=npa/d^2$. The limit of this expression, as $t\to 1$,  always
exists. In particular $\#\{i:\ \chi(K_2^{(i)})\not=1\}\geq 2$ (cf. also
with \ref{g3}(1)). If this number is strict greater than 2, then the limit
is zero. If $\chi(K_2^{(i)})\not=1$ exactly for two indices $i$ and $j$, then
using \ref{g3}(1) clearly $\chi(K_2^{(i)})\chi(K_2^{(j)})=1$. Since there
are exactly $d(d-1)/2$ such pairs, one gets:
$$\et_{M_2,\sigma_{can}}(1)=\frac{1}{|H_1(M_2)|}\cdot
\lim_{t\to 0}\frac{(t^\alpha-1)^d}{(t^{d\alpha/n}-1)
(t^{d\alpha/p}-1)(t^{\alpha/a}-1)^{d-2}}
\cdot \sum_{i\not=j}\sum_{\xi\in \Z_a^*}\frac{1}{(\xi-1)(\bar{\xi}-1)}$$
\begin{equation*}
=\frac{np}{ad^2}\sum_{i\not=j}\sum_{\xi\in \Z_a^*}\frac{1}{(\xi-1)(\bar{\xi}-1)}
=\frac{np}{ad^2}\frac{d(d-1)}{2}
\sum_{\xi\in \Z_a^*}\frac{1}{(\xi-1)(\bar{\xi}-1)}=
\frac{np(d-1)(a^2-1)}{24ad},
\tag{$*$}
\end{equation*}
since
\begin{equation*}
\sum_{\xi\in \Z_a^*}\frac{1}{(\xi-1)(\bar{\xi}-1)}=\frac{a^2-1}{12}.
\tag{$**$}
\end{equation*}
Consider now the manifold $M$. Then using \ref{s7}(B) (and/or its proof),
by the same argument as  above, one obtains:
\begin{equation*}
\et_{M,\sigma_{can}}(1)=%\sum_i\et_{M_2^{(i)},\sigma_{can}}(1)+
\frac{np}{ad^2}\sum_{i\not=j}\sum_{\xi\in \Z_a^*}
\frac{\Delta_{M_1^{(i)}}(\xi)\cdot \Delta_{M_1^{(j)}}(\bar{\xi})}{(\xi-1)
(\bar{\xi}-1)}.
\tag{$***$}
\end{equation*}
Finally, making the difference between ($***$) and ($*$) one gets the result.
\end{proof}

\subsection{Example/Discussion.}\label{g9} Assume
that $M$ satisfies {\em WA1'}
and {\em WA2'}, and additionally $(M_1^{(i)},K_1^{(i)})=(M_1,K_1)$
for some {\em integral homology sphere} $M_1$. Then
$${\mathcal O}(\lambda_W)/2= \frac{np(d-1)}{2ad}\cdot
(\Delta^{\no}_{M_1})''(1);$$
$${\mathcal O}(\et_{\cdot ,\sigma_{can}}(1))=
\frac{np(d-1)}{2ad}\cdot \sum_{\xi\in \Z_a^*}\
\frac{\Delta_{M_1}(\xi)\cdot \Delta_{M_1}(\bar{\xi})-1}
{(\xi-1)(\bar{\xi}-1)}.$$
Recall that the modified Seiberg-Witten invariants $\ssw^0_M(\sigma_{can})$
is defined by the difference $\et_{M,\sigma_{can}}(1)-\lambda_W(M)/2$
(cf. \ref{sw}(13)).
Notice the remarkable fact that in the above expressions the coefficients
  before the Alexander invariants became the same.  Hence
\begin{equation*}
{\mathcal O}(\ssw^0_{\cdot}(\sigma_{can}))=\frac{np(d-1)}{2ad}{\mathcal D}_a,\
\ \mbox{where}  \ \ {\mathcal D}_a:= \sum_{\xi\in \Z_a^*}\
\frac{\Delta_{M_1}(\xi)\cdot \Delta_{M_1}(\bar{\xi})-1}
{(\xi-1)(\bar{\xi}-1)} - (\Delta^{\no}_{M_1})''(1).
\tag{${\mathcal D}$}
\end{equation*}
Recall that $\Delta^{\no}_{M_1}(t)$ is a symmetric polynomial
(cf. \cite[2.3.1]{Lescop}) with $\Delta^{\no}_{M_1}(1)=1$ (cf. \ref{s0e}(6)).
In the sequel we will compute
explicitly ${\mathcal D}_a$, provided that $a$ is sufficiently large,
in terms of the coefficients $\{c_i\}_{i=1}^r$ of $\Delta^{\no}_{M_1}(t)$
(cf. \ref{def1}(3)).

The contribution $(\Delta^{\no}_{M_1})''(1)$ is easy: it is
$ \sum_{i=1}^r\ 2i^2\, c_i$.  By \ref{s0e}(9),
$\det(V)=1$, hence by \ref{s0e}(7),
$\Delta_{M_1}(t)=t^r\cdot  \Delta^{\no}_{M_1}(t)$.
In particular, $\Delta_{M_1}(\xi)\cdot \Delta_{M_1}(\bar{\xi})=
  \Delta^{\no}_{M_1}(\xi)\cdot   \Delta^{\no}_{M_1}(\bar{\xi})$.
Then write
$$\frac{\Delta^{\no}_{M_1}(t)}{1-t}=\frac{1}{1-t}-\sum_{i=1}^rc_i(1+t+\cdots+
t^{i-1}) +\sum_{i=1}^rc_i(t^{-1}+\cdots t^{-i}).$$
Using the identity $\sum_{\xi\in \Z_a^*}1/(1-\xi)=(a-1)/2$, an elementary
computation gives
$$\sum_{\xi\in \Z_a^*}\
\frac{\Delta^{\no}_{M_1}(\xi)\cdot \Delta^{\no}_{M_1}(\bar{\xi})-1}
{(1-\xi)(1-\bar{\xi})}= \sum_{i=1}^r\ 2i^2\, c_i +2a\cdot
{\mathcal D}(\Delta^{\no}_{M_1}(t)), \ \mbox{provided that $a\geq 2r$}.$$
In particular, if $a\geq 2r$, then
${\mathcal D}_a=2a\cdot {\mathcal D}(\Delta^{\no}_{M_1}(t))$. Hence
%2a\cdot \Big[\sum_{i,j=1}^r \, c_ic_j\min(i,j)-\sum_{i=1}^r\, i\, c_i\Big],$$
$${\mathcal O}(\ssw^0_{\cdot}(\sigma_{can}))=\frac{np(d-1)}{d}\cdot
{\mathcal D}(\Delta^{\no}_{M_1}(t)).$$
This raises the following natural question: for what Alexander polynomials
the expression ${\mathcal D}(\Delta^{\no}_{M_1}(t))$ is zero?
The next lemma provides such an example (the proof is elementary
and it is left to the reader).

\subsection{Lemma.}\label{g11} {\em If $\{c_i\}_{i=1}^r$ is an alternating set
then ${\mathcal D}(c_1,\ldots, c_r)=0$.}

\vspace{2mm}

The above discussions have the following topological consequence:

 Fix two relative prime positive integers $p$ and $a$.
Let $K$ be the primitive simple curve in $\partial T(L_1)$ with homology class
$am_1+p\lambda_1$.  Let $M$ denote
the $n$-cyclic cover of $N_1$ branched along $K$.
Set $d:=\mbox{gcd}(n,p)$, and let $M_1$ be the $(n/d)$-cyclic cover of
$N_1$ branched along $L_1$. Denote by $K_1$ the preimage of $L_1$ via
this cover. Finally,
let $\Delta^{\no}_{M_1}(t)$ be the  normalized Alexander polynomial of
$(M_1,K_1)$.
\subsection{Corollary.}\label{swc} {\em Consider the above data. Additionally,
assume that WA2' is satisfied and $M$ is a rational homology sphere. Then:

(A) $(d-1)\cdot (\mbox{gcd}(n,a)-1)=0$.

(B) If $d=1$, then
$$\ssw^0_M(\sigma_{can})=\ssw^0_{M_1}(\sigma_{can})+
\ssw^0_{\Sigma(p,a,n)}(\sigma_{can}).$$

(C) If gcd$(n,a)=1$, $a\geq \deg\Delta_{M_1}(t)$, and
$M_1$ is an integral homology sphere,  then
$$\ssw^0_M(\sigma_{can})=d\cdot \ssw^0_{M_1}(\sigma_{can})+
\ssw^0_{\Sigma(p,a,n)}(\sigma_{can})+\frac{np(d-1)}{d}\cdot {\mathcal D}
(\Delta^{\no}_{M_1}(t)).$$
If the coefficients of $\Delta^{\no}_{M_1}(t)$ form an alternating set,
then ${\mathcal D}(\Delta^{\no}_{M_1}(t))=0$. }
\begin{proof}
Consider the following schematic splicing of splice diagrams (cf. \ref{spl}):

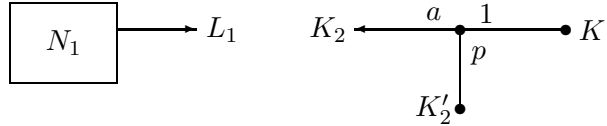
\begin{picture}(300,60)(-50,10)
\put(20,30){\framebox(40,30){$N_1$}}
\put(60,50){\vector(1,0){30}}
\put(190,50){\vector(-1,0){40}}
\put(100,50){\makebox(0,0){$L_1$}}
\put(140,50){\makebox(0,0){$K_2$}}
\put(240,50){\makebox(0,0){$K$}}
\put(190,50){\circle*{4}}
\put(190,20){\circle*{4}}
\put(230,50){\circle*{4}}
\put(190,50){\line(1,0){40}}
\put(190,50){\line(0,-1){30}}
\put(180,20){\makebox(0,0){$K_2'$}}
\put(180,56){\makebox(0,0){$a$}}
\put(200,56){\makebox(0,0){$1$}}
\put(197,40){\makebox(0,0){$p$}}
\end{picture}

The result of the splicing can be identified with $N_1$ (and under this
identification $L_1$ is identified with $K_2'$). The advantage of this
splicing representation is that it emphasizes the position of the
knot $K$ in the Seifert component $\Sigma(p,a,1)$. If $M$ is a rational
 homology sphere then the $n$-cyclic cover of $\Sigma(p,a,1)$
branched along $K$ (which is $\Sigma(p,a,n)$) should be rational homology
sphere, hence (A) follows. If $d=1$ then $M$ has a splice decomposition of
the following schematic plumbing  diagrams
 (where at the right $M_2=\Sigma(p,a,n)$
and the dots mean gcd$(n,a)$ arms):

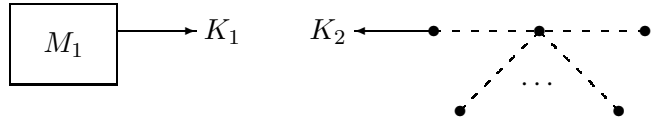
\begin{picture}(300,55)(-50,15)
\put(20,30){\framebox(40,30){$M_1$}}
\put(60,50){\vector(1,0){30}}
\put(180,50){\vector(-1,0){30}}
\put(100,50){\makebox(0,0){$K_1$}}
\put(140,50){\makebox(0,0){$K_2$}}
\put(220,50){\circle*{4}}
\put(180,50){\circle*{4}}
\put(250,20){\circle*{4}}
\put(190,20){\circle*{4}}
\put(260,50){\circle*{4}}
\dashline{3}(180,50)(260,50)
\dashline{3}(220,50)(250,20)
\dashline{3}(220,50)(190,20)
\put(220,30){\makebox(0,0){$\cdots$}}
\end{picture}

%\dashline{3}(50,60)(175,60)

\noindent Here $o_1=o_2=1$ and $k_1=k_2=0$, and $A$ is the identification
$\lambda_2=m_1$, $m_2=\lambda_1$. Therefore, part (B) follows from
\ref{s2}(4) and \ref{s7}(A). The last case corresponds exactly to the
situation treated in \ref{g9}.
\end{proof}

\subsection{Remarks.}\label{g12}  (1) Our final goal
(see  the following sections)  is to prove the additivity result
${\mathcal O}(\ssw^0(\sigma_{can}))=0$ for any $(M_1,K_1)$,
which can be represented as a cyclic
cover  of $S^3$ branched along the  link $K_f\subset S^3$ of an
arbitrary irreducible
(complex) plane curve singularity
(even if $M_1$ is not an integral    homology sphere),
provided that $a$ is sufficiently large.
This means that from the above Corollary, part (C), we will need to eliminate
the assumption about the vanishing of $H_1(M_1)$. The assumption
about $a$ will follow from the  special property  \ref{a1}(6)) of
irreducible plane curve singularities (cf. also with \ref{a4}, especially
with (7)).

The proof of the vanishing of  the ${\mathcal D}$-correction  will take
up most of the last section of the paper. It relies in a crucial manner on
the alternating nature of the Alexander polynomial $\Delta^{\no}_{S^3}(K_f)(t)$
 of any irreducible plane curve singularity $f$, fact
which will be establish in  Proposition \ref{a2}.

(2) It is really interesting and remarkable, that the
behavior of the modified Seiberg-Witten  invariant with respect to (some)
splicing and cyclic covers (constructions, which basically are
 topological in nature)
definitely gives  preference to the Alexander polynomials of some algebraic
links. The authors hope that a better understanding of this phenomenon
would lead to some deep properties of the Seiberg-Witten invariant.

\subsection{Example.}\label{dnotzero} \ In general, in \ref{swc},
the invariant ${\mathcal D}
(\Delta^{\no}_{M_1}(K_1)(t))$ does not vanish. In order to see this, start for
example with a pair $(N_1,L_1)$ with non-zero ${\mathcal D}
(\Delta^{\no}_{N_1}(L_1)(t))$, and consider the case when $d|n$.
(If $d\not=1$, then the coefficient of ${\mathcal D}
(\Delta^{\no}_{N_1}(L_1)(t))$ in \ref{swc}(C) will  be non-zero as well.)

Next, we show how one can construct a pair $(N,L)$ which satisfies WA2,
$H_1(N)=0$, but ${\mathcal D}(\Delta^{\no}_{N}(L)(t))\not=0$.
First, we notice the following fact.

If the Alexander polynomial $\Delta^{\no}_M(K)(t)$ is realizable for some
pair $(M,K)$ (satisfying WA2 and $H_1(M)=0$), then the $k$-power of this
polynomial is also realizable for some pair $(M^k,K^k)$ (satisfying WA2 and
$H_1(M^k)=0$). Indeed, assume that $(M,K)$ has a schematic splice diagram
of the following form:

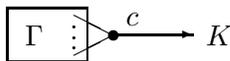
\begin{picture}(400,35)(-100,37)
\put(30,40){\framebox(30,20)}
\put(70,50){\circle*{4}}
\put(70,50){\line(-2,1){15}}
\put(70,50){\line(-2,-1){15}}
\put(55,52){\makebox(0,0){$\vdots$}}
\put(70,50){\vector(1,0){30}}
\put(77,56){\makebox(0,0){$c$}}
\put(110,50){\makebox(0,0){$K$}}
\put(40,50){\makebox(0,0){$\Gamma$}}
\end{picture}

\noindent Then let $(M^k,K^k)$  be given by the following schematic splice
diagram:

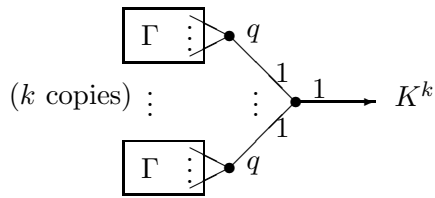
\begin{picture}(400,85)(-100,37)
\put(30,40){\framebox(30,20)}
\put(70,50){\circle*{4}}
\put(70,50){\line(-2,1){15}}
\put(70,50){\line(-2,-1){15}}
\put(55,52){\makebox(0,0){$\vdots$}}
\put(70,50){\line(1,1){25}}
\put(79,51){\makebox(0,0){$q$}}
\put(40,50){\makebox(0,0){$\Gamma$}}

\put(30,90){\framebox(30,20)}
\put(70,100){\circle*{4}}
\put(70,100){\line(-2,1){15}}
\put(70,100){\line(-2,-1){15}}
\put(55,102){\makebox(0,0){$\vdots$}}
\put(70,100){\line(1,-1){25}}
\put(79,101){\makebox(0,0){$q$}}
\put(40,100){\makebox(0,0){$\Gamma$}}

\put(40,77){\makebox(0,0){$\vdots$}}
\put(80,77){\makebox(0,0){$\vdots$}}
\put(10,77){\makebox(0,0){($k$ copies)}}
\put(95,75){\circle*{4}}
\put(95,75){\vector(1,0){30}}
\put(104,80){\makebox(0,0){$1$}}
\put(90,86){\makebox(0,0){$1$}}
\put(90,65){\makebox(0,0){$1$}}
\put(140,78){\makebox(0,0){$K^k$}}
\end{picture}

\noindent Here, we take $q$ sufficiently large (in order to assure that
the new edges will also satisfy the  algebraicity condition \cite{EN}(9.4)),
and also $q$ should be relative prime with some integers which appear
as  decorations of $\Gamma$ (see [loc. cit.]). Obviously, by construction,
$M^k$ is an integral homology sphere. Then, by [loc. cit.] 12.1,
one can easily verify that
$$\Delta^{\no}_{M^k}(K^k)(t)= \Delta^{\no}_{M}(K)(t)^k.$$
For example, if $(M,K)=(S^3,K_f)$, where $K_f$ is the $(2,3)$-torus knot
(or, equivalently, the knot of the plane curve singularity $f=x^2+y^3$,
cf. \ref{a1}),
then $\Delta^{\no}_{S^3}(K_f)(t)=t-1+1/t$ (see \ref{a1}(5)).
Now, if we take $k=2$ and $q=7$ then
$(M^2,K^2)$ has the following splice, respectively plumbing graph:

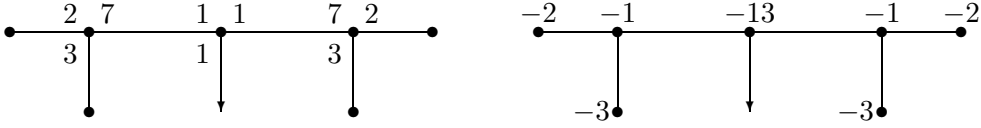
\begin{picture}(400,60)(0,0)
\put(20,40){\circle*{4}}
\put(50,40){\circle*{4}}
\put(100,40){\circle*{4}}
\put(150,40){\circle*{4}}
\put(180,40){\circle*{4}}

\put(50,10){\circle*{4}}
\put(150,10){\circle*{4}}
\put(100,40){\vector(0,-1){30}}
\put(50,40){\line(0,-1){30}}
\put(20,40){\line(1,0){160}}
\put(150,40){\line(0,-1){30}}

\put(43,47){\makebox(0,0){$2$}}
\put(57,47){\makebox(0,0){$7$}}
\put(93,47){\makebox(0,0){$1$}}
\put(107,47){\makebox(0,0){$1$}}
\put(143,47){\makebox(0,0){$7$}}
\put(157,47){\makebox(0,0){$2$}}
\put(43,32){\makebox(0,0){$3$}}
\put(93,32){\makebox(0,0){$1$}}
\put(143,32){\makebox(0,0){$3$}}

\put(220,40){\circle*{4}}
\put(250,40){\circle*{4}}
\put(300,40){\circle*{4}}
\put(350,40){\circle*{4}}
\put(380,40){\circle*{4}}

\put(250,10){\circle*{4}}
\put(350,10){\circle*{4}}
\put(300,40){\vector(0,-1){30}}
\put(250,40){\line(0,-1){30}}
\put(220,40){\line(1,0){160}}
\put(350,40){\line(0,-1){30}}

\put(220,47){\makebox(0,0){$-2$}}
\put(250,47){\makebox(0,0){$-1$}}
\put(300,47){\makebox(0,0){$-13$}}
\put(350,47){\makebox(0,0){$-1$}}
\put(380,47){\makebox(0,0){$-2$}}
\put(240,10){\makebox(0,0){$-3$}}
\put(340,10){\makebox(0,0){$-3$}}

\end{picture}

\noindent Then $(M^2,K^2)$ is algebraic, $H_1(M^2)=0$. But
$\Delta^{\no}_{M^2}(K^2)(t)=(t-1+1/t)^2$ whose coefficients are not
alternating. In fact, $r=2$, $c_1=-2$ and $c_2=1$; in particular
${\mathcal D}(\Delta^{\no}_{M^2}(K^2)(t))=2$.

\vspace{2mm}

We end this section with the following property
 which is needed in the last section.

\subsection{Lemma.}\label{detg} {\em Assume that $M$ satisfies
WA1' and WA2'
 with $(M_1^{(i)},K_1^{(i)})=(M_1,K_1)$. Let $\Gamma$ denote
the plumbing graph of $M$. Let $v$ be the central vertex of $M_2$ considered in
$M$, and let $\Gamma_-$ be that connected component of $\Gamma\setminus \{v\}$
which contains the vertices of $M_1^{(1)}$. Then $|\det(\Gamma_-)|=
a\cdot |H_1(M_1)|$.}

\begin{proof}
If $I$ denotes the intersection matrix of $M$, then
\[
-I^{-1}_{vv}\stackrel{\ref{s0d}(5)}{=}Lk_M(g_v,g_v)\stackrel{\ref{s5}(12)}{=} Lk_{\Sigma(p,a,n)}(O,O)\stackrel{\ref{g3} (L\#c)}{=}npa/d^2.
\]
On the other
hand, $I^{-1}_{vv}$ can be computed from the determinants of the components of
$\Gamma\setminus \{v\}$, hence (cf. also with \ref{g2}(a))
\[
-|H_1(M)|\cdot I^{-1}_{vv}=|\det(\Gamma_-)|^d\cdot pn/d^2 .
\]
By \ref{s5}(7) and \ref{g2}(e)  $|H_1(M)|=|H_1(M_1)|^d\cdot a^{d-1}$, hence
the result follows.
\end{proof}

\section{Properties of irreducible plane curve singularities}

\subsection{The topology of an irreducible plane curve singularity.}\label{a1}
Consider an irreducible plane   curve singularity $f:(\C^2,0)\to (\C,0)$
with Newton pairs $\{(p_k,q_k)\}_{k=1}^s$ (cf. \cite{EN}, page 49).
Clearly gcd$(p_k,q_k)=1$ and $p_k\geq 2$ and $q_k\geq 2$.
Define the integers $\{a_k\}_{k=1}^s$ by
\begin{equation*}
a_1=q_1 \ \mbox{and } \  a_{k+1}=q_{k+1}+p_{k+1}p_ka_k\ \ \mbox{if $k\geq 1$}.
\tag{1}
\end{equation*}
Then again,  gcd$(p_k,a_k)=1$ for any $k$.
The minimal (good) embedded resolution graph of the pair $(\C^2,\{f=0\})$
has the following schematic form:

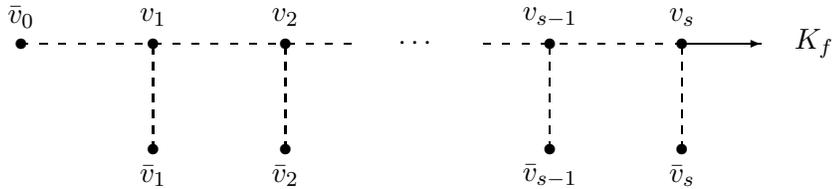
\begin{picture}(400,80)(0,0)
\put(50,60){\circle*{4}}
\put(100,60){\circle*{4}}
\put(150,60){\circle*{4}}
\put(250,60){\circle*{4}}
\put(300,60){\circle*{4}}
\put(100,20){\circle*{4}}
\put(150,20){\circle*{4}}
\put(250,20){\circle*{4}}
\put(300,20){\circle*{4}}
\put(50,70){\makebox(0,0){$\bar{v}_0$}}
\put(100,70){\makebox(0,0){$v_1$}}
\put(150,70){\makebox(0,0){$v_2$}}
\put(250,70){\makebox(0,0){$v_{s-1}$}}
\put(300,70){\makebox(0,0){$v_s$}}
\put(100,10){\makebox(0,0){$\bar{v}_1$}}
\put(150,10){\makebox(0,0){$\bar{v}_2$}}
\put(250,10){\makebox(0,0){$\bar{v}_{s-1}$}}
\put(300,10){\makebox(0,0){$\bar{v}_s$}}
\put(350,60){\makebox(0,0){$K_f$}}
\dashline{3}(50,60)(175,60)
\dashline{3}(225,60)(300,60)
\put(100,20){\dashbox{3}(0,40){}}
\put(150,20){\dashbox{3}(0,40){}}
\put(250,20){\dashbox{3}(0,40){}}
\put(300,20){\dashbox{3}(0,40){}}
\put(200,60){\makebox(0,0){$\cdots$}}
\put(300,60){\vector(1,0){30}}
\end{picture}

\noindent This can be identified with the plumbing graph $\Gamma(S^3,K_f)$,
where $K_f$ is the link of $f$ (with only one component)
in the Milnor sphere $S^3$. In the above diagram
we emphasized only  those vertices $\{\bar{v}_k\}_{k=0}^s$
and $\{v_k\}_{k=1}^s$ which have  $\bar{\delta}\not=2$.
We denote the  set of these vertices by ${\mathcal V}^*$.
The dash-line between two such vertices replaces a string
\begin{picture}(100,10)(0,3)
\put(20,5){\circle*{4}}
\put(40,5){\circle*{4}}
\put(80,5){\circle*{4}}
\put(10,5){\line(1,0){40}}
\put(70,5){\line(1,0){20}}
\put(62,5){\makebox(0,0){$\cdots$}}
\end{picture}
. In our discussion the corresponding self-intersection (or Euler) numbers
will be not important (the interested reader can find the complete
description of the graph in \cite{EN} section 22, or in \cite{NDedII}).
The above numerical data $\{(p_k,a_k)\}_k$
and the set of vertices ${\mathcal V}^*$ is codified in the
splice  diagram (cf. \cite{EN}):

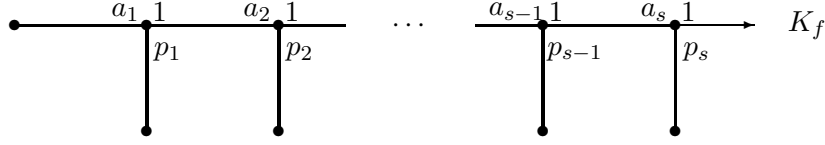
\begin{picture}(400,65)(0,10)
\put(50,60){\circle*{4}}
\put(100,60){\circle*{4}}
\put(150,60){\circle*{4}}
\put(250,60){\circle*{4}}
\put(300,60){\circle*{4}}
\put(100,20){\circle*{4}}
\put(150,20){\circle*{4}}
\put(250,20){\circle*{4}}
\put(300,20){\circle*{4}}
\put(350,60){\makebox(0,0){$K_f$}}
\put(92,65){\makebox(0,0){$a_1$}}
\put(142,65){\makebox(0,0){$a_2$}}
\put(240,65){\makebox(0,0){$a_{s-1}$}}
\put(292,65){\makebox(0,0){$a_s$}}

\put(105,65){\makebox(0,0){$1$}}
\put(155,65){\makebox(0,0){$1$}}
\put(255,65){\makebox(0,0){$1$}}
\put(305,65){\makebox(0,0){$1$}}

\put(108,50){\makebox(0,0){$p_1$}}
\put(158,50){\makebox(0,0){$p_2$}}
\put(262,50){\makebox(0,0){$p_{s-1}$}}
\put(308,50){\makebox(0,0){$p_s$}}
\put(200,60){\makebox(0,0){$\cdots$}}
\put(50,60){\framebox(125,0){}}
\put(225,60){\framebox(75,0){}}
\put(100,20){\framebox(0,40){}}
\put(150,20){\framebox(0,40){}}
\put(250,20){\framebox(0,40){}}
\put(300,20){\framebox(0,40){}}
\put(300,60){\vector(1,0){30}}
\end{picture}

\noindent The knot $K_f\subset S^3$ defines a set of weights $\{w_{v}(u)\}
_{v\in {\mathcal V}^*}$
as in \ref{s0d} (where $K_f=g_u$, and evidently $o(u)=1$). In terms of
the resolution,
$w_v(u)$ is exactly the vanishing order (multiplicity) of $f\circ\pi$
along the exceptional divisor codified by $v$, where $\pi$ denoted the
resolution map. In the sequel we write just $w_v$ for it. Then (cf.
\cite{EN}, section 10) one has:
\begin{equation*}
\begin{array}{ll}
w_{v_k}=a_kp_kp_{k+1}\cdots p_s & \mbox{for $1\leq k\leq s$};\\
w_{\bar{v}_0}=p_1p_2\cdots p_s; & \\
w_{\bar{v}_k}=a_kp_{k+1}\cdots p_s & \mbox{for $1\leq k\leq s$}. \end{array}
\tag{2}
\end{equation*}
Recall that the characteristic polynomial
$\Delta(f)(t):=\Delta_{S^3}(K_f)(t)$ of the
monodromy acting on the first homology of the Milnor fiber of $f$
is given by A'Campo's formula \ref{s0e}(8):
\begin{equation*}
\frac{\Delta_{S^3}(K_f)(t)}{t-1}= \prod_{v\in{\mathcal V^*}}\
(t^{w_v}-1)^{\bar{\delta}_v-2}.
\tag{3}
\end{equation*}
In inductive proofs and constructions
(over the number of Newton pairs of $f$),
it is convenient to use the notation $f_{(l)}$ for an irreducible
plane curve singularity with Newton pairs $\{(p_k,q_k)\}_{k=1}^l$,
where $1\leq l\leq s$. Evidently $f_{(s)}=f$, and $f_{(1)}$ can be taken
as the Brieskorn singularity $x^{p_1}+y^{a_1}$.
We write $\Delta(f_{(l)}) $  for the characteristic polynomial associated
with $f_{(l)}$. Then from (2) and (3) one gets
\begin{equation*}
\Delta(f_{(l)})(t)=\Delta(x^{p_l}+y^{a_l})(t)\cdot \Delta(f_{(l-1)})(t^{p_l})
\ \ \mbox{for $l\geq 2$},
\tag{4}
\end{equation*}
where
\begin{equation*}
\Delta(x^{p}+y^{a})(t)=\frac{(t^{pa}-1)(t-1)}{(t^p-1)(t^a-1)}.
\tag{5}
\end{equation*}
By induction, using the identities (1), one can prove
(see e.g. \cite[5.2]{NDedI})
\begin{equation*}
a_l>p_l\cdot \deg\Delta(f_{(l-1)})\ \ \mbox{for any $l\geq 2$}.
\tag{6}
\end{equation*}

\subsection{Proposition.}\label{a2} {\em
(a) $\Delta(f)(0)=\Delta(f)(1)=1$, and the degree of $\Delta(f)(t)$
is even (say $2r$).

(b) If $ \Delta(f)(t)=\sum_{i=0}^{2r}b_it^i$, then the set $\{b_i\}_{i=0}^{2r}$
is alternating (cf. \ref{def1}(4)).

(c) The coefficients $\{c_i\}_{i=1}^r$ of $\Delta^{\no}(f)(t):=t^{-r}
\Delta(f)(t)$ (cf. \ref{def1}(3)) are alternating as well.}

\begin{proof} (a) is clear from (4) and (5), and (c)  follows easily  from (b).
We will prove (b) by induction over $s$.
For each $1\leq l\leq s$ we verify that there exist

(i) $a_l$-residue classes $\{r_1,\ldots, r_t\}\subset \{1,2, \ldots, a_l-1\}$
(where $t$ may depend on $l$); and

(ii) \ integers $n_1,\ldots, n_t\in \N$, such that
$$\Delta(f_{(l)})(t)=1+\sum_{i=1}^t\sum_{j=0}^{n_i}\, t^{r_i+ja_l}\cdot (t-1).$$
It is clear that the coefficients of a polynomial of this form are alternating.

Let us start with the case $l=1$. Write $(p_1,a_1)=(p,a)$. Then (cf.  (5))

\[
\Delta(x^p+y^a)(t)=(t^{p(a-1)}+\cdots +t^p+1)/Q(t),\;\;Q(t):=t^{a-1}+\cdots +t+1.
\]
 For each $i=0, 1,\ldots,a-1$ write $pi$ in the form
$x_ia+r_i$ for some $r_i\in \{0, \ldots, a-1\}$. Since gcd$(p,a)=1$,
$\{r_i\}_i=\{0,\ldots,a-1\}$, and $p|r_i$ if and only if $x_i=0$.
Therefore,
$$\sum _{i=0}^{a-1}t^{pi}=
Q(t)+ \sum_{i:\,p\,\nmid \,r_i}t^{r_i}(t^{x_ia}-1)=Q(t)\cdot
\Big[ 1+\sum_{i:\,p\,\nmid\, r_i}\sum_{j=0}^{x_i-1}\, t^{r_i+ja}(t-1)\Big].$$
Now we prove that $\Delta(f_{(l)}(t)$ has  a similar form. By the inductive
step, assume that
$$\Delta(f_{(l-1)})(t)=1+\sum_{i=1}^t\sum_{j=0}^{n_i}\, t^{r_i+ja_{l-1}}
(t-1).$$
Then, using \ref{a1}(4) and (5), one gets for $\Delta(f_{(l)})(t)$:
$$\frac{(t^{p_la_l}-1)(t-1)}{(t^{p_l}-1)(t^{a_l}-1)}+
\sum_{i=1}^t\sum_{j=0}^{n_i}\, t^{(r_i+ja_{l-1})p_l}
\cdot \frac{(t^{p_la_l}-1)(t-1)}{t^{a_l}-1}.$$
Let $\{s_j\}_{j=0}^{a_l-1}$ be the set of $a_l$-residues classes. Then
(using the result of case $l=1$) the above expression  reads as
$$1+\sum_{j:\,p_l\,\nmid\, s_j}\sum_{k=0}^{x_j-1}\, t^{s_j+ka_l}(t-1)+
\sum_{i=1}^t\sum_{j=0}^{n_i}\sum_{k=0}^{p_l-1}\, t^{(r_i+ja_{l-1})p_l+ka_l}
\cdot (t-1).$$
Notice that \ref{a1}(6) guarantees that for each $i$ and $j$ one has the
inequality $(r_i+ja_{l-1})p_l<a_l$, hence these numbers can be considered as
(non-zero)
$a_l$-residues classes. Moreover, they are all different from the
residue classes  $\{s_j: p_l\,\nmid\, s_j\}$ since they are all divisible by
$p_l$.
\end{proof}

The ``alternating property'' of the coefficients of
the Alexander polynomial of any  irreducible
plane curve singularity will be crucial in the computation of the
Reidemeister-Turaev torsion of $\{f+z^n=0\}$.
 The key algebraic fact is summarized in the next property:

\subsection{Algebraic Lemma.}\label{a3} {\em In the next expressions
  $t$ is a free variable  and $a$ is a positive
integer. $\Z_a$ is identified with the $a$-roots of unity.
 Assume that the coefficients of a polynomial $\Delta(t)\in \Z[t]$
form an alternating set, $\Delta(1)=1$, and  $a\geq \deg\Delta$. Then:

(a) For an arbitrary complex number $A$ one has:
$$\frac{1}{a}\sum_{\xi\in \Z_a}\frac{\Delta(\xi t)}{1-\xi t}\cdot
\frac{\Delta(\bar{\xi}At)}{1-\bar{\xi}At}=
\frac{(1-A^at^{2a})\cdot \Delta(At^2)}{(1-t^a)(1-A^at^a)(1-At^2)}.$$

(b)  For arbitrary integers $d\geq 2$ and $k\geq 1$ one has:
$$\frac{1}{a^{d-1}}\sum_{\xi_1,\ldots, \xi_d\in \Z_a
\atop \xi_1\cdots\xi_d=1}\frac{\Delta(\xi_1t)}{1-\xi_1 t}
\cdot\, \cdots\, \cdot
\frac{\Delta(\xi_{d-1}t)}{1-\xi_{d-1}t}
\cdot \frac{\Delta(\xi_dt^k)}{1-\xi_dt^k}=
\frac{(1-t^{a(d+k-1)})\cdot \Delta(t^{d+k-1})}{(1-t^a)^{d-1}(1-t^{ak})
(1-t^{d+k-1})}.$$}

\begin{proof} The assumption about $\Delta(t)$ guarantees that one can write
$\Delta(t)=1-R(t)(1-t)$ for some $R(t)=\sum_{j\geq1}\tilde{b}_jt^j$
with $\tilde{b}_j\in \{0,1\}$ for all $j$. Then the left hand side of (a) is
$$\frac{1}{a}\sum_{\xi\in \Z_a}\frac{1}{(1-\xi t)(1-\bar{\xi}At)}-
\frac{1}{a}\sum_{\xi\in \Z_a}\frac{R(\bar{\xi}At)}{1-\xi t}-
\frac{1}{a}\sum_{\xi\in \Z_a}\frac{R(\xi t)}{1-\bar{\xi}A t}+
\frac{1}{a}\sum_{\xi\in \Z_a}R(\xi t)\cdot R(\bar{\xi}At).$$
The first sum (with the coefficient $1/a$) can be written in the form
$$\frac{1}{a}(1+\xi t+\xi^2t^2+\cdots)(1+\bar{\xi}At+\bar{\xi}^2A^2t^2+\cdots)
=\frac{1}{a}\sum_{n\geq 0}\, \sum_{j=0}^n\, \xi^{n-2j}A^jt^n.$$
This, by an elementary computation gives:
$$\frac{1-A^at^{2a}}{(1-t^a)(1-A^at^a)(1-At^2)}.$$
The second term gives $R(At^2)/(1-t^a)$. In order to prove this, first
notice that the formula is additive in the polynomial $R$, hence it is enough
to verify for $R(t)=t^k$ for all $0\leq k< a$. The case $k=0$ is easy,
it is equivalent with the identity
\begin{equation*}
\frac{1}{a}\sum_{\xi}\frac{1}{1-\xi t}=\frac{1}{1-t^a}.
\tag{$*$}
\end{equation*}
If $1\leq k<a$, then write
$$\frac{1}{a}\sum_{\xi}\frac{\bar{\xi}^kA^kt^k}{1-\xi t}=
\frac{A^kt^{2k}}{a}\Big[\sum_{\xi}\frac{1}{1-\xi t}+\sum _{\xi}
(\xi t)^{-1}+\cdots +(\xi t)^{-k}\Big].$$
Since $k<a$ the last sum is zero  (here $k<a$ is crucial !),
hence ($*$) gives the claimed identity.

By similar method,
the third term is $R(At^2)/(1-A^at^a)$. Finally,  the forth is $R(At^2)$
(here one needs to apply the alternating property, namely that
$\tilde{b}_j^2=\tilde{b}_j$ for any $j$).

For part (b), use (a) and induction over $d$.
For this, write $\xi_d$ as $\bar{\xi}_1\cdots \bar{\xi}_{d-1}$ and use (a) for
$\xi=\xi_{d-1}$ and $A=\bar{\xi}_1\cdots\bar{\xi}_{d-2}t^{k-1}$. Then apply
the inductive step.
\end{proof}

%(b) For any integer $d\geq 2$, and integer $k\geq 1$  one has:
%$$\frac{1}{a^{d-1}}\sum_{\xi_1,\ldots, \xi_d\in \Z_a
%\atop \xi_1\cdots\xi_d=1}\frac{1}{1-\xi_1 t}\cdot\, \cdots\, \cdot
%\frac{1}{1-\xi_{d-1}t}\cdot \frac{1}{1-\xi_dt^k}=
%\frac{1-t^{a(d+k-1)}}{(1-t^a)^{d-1}(1-t^{ak})(1-t^{d+k-1})}.$$

\subsection{Remarks.}\label{a4}

\noindent (1)  Let $\Delta(t)$ and $a$ be as in \ref{a3}.
 The expression in \ref{a3}(a), with $A=1$,
 has a pole of order 2 at $t=1$. This comes from the pole of the summand
given by $\xi=1$.  Therefore:
\[\frac{1}{a}\sum_{\xi\in \Z_a^*}\frac{\Delta(\xi)}{1-\xi}\cdot
\frac{\Delta(\bar{\xi})}{1-\bar{\xi}}=
\lim_{t\to 1} \Big[\,
\frac{(1-t^{2a})\cdot \Delta(t^2)}{(1-t^a)^2(1-t^2)}-
\frac{\Delta(t)^2}{a(1-t)^2}\, \Big]
\]
\[
=\frac{a^2-1}{12a}+\frac{1}{a}\Big[ \Delta'(1)-\Delta'(1)^2+\Delta''(1)\Big],
\]
where the first equality follows from \ref{a3}, the second by a  computation.

(2) Assume that $\Delta(t)$ is an arbitrary symmetric polynomial
of degree $2r$, and write
$\Delta^{\no}(t)=t^{-r}\Delta(t)$. Then it is easy to show that
$$\Delta'(1)=r\Delta(1), \ (\Delta^{\no})'(1)=0, \ \mbox{and}
\ \
(\Delta^{\no})''(1)= (r-r^2) +\Delta''(1)/\Delta(1).$$

(3) If one combines (1) and (2), then for a symmetric polynomial $\Delta(t)$
with alternating coefficients and with $\Delta(1)=1$ one gets
$$\sum_{\xi\in \Z_a^*}\frac{\Delta(\xi)}{1-\xi}\cdot
\frac{\Delta(\bar{\xi})}{1-\bar{\xi}}=
\frac{a^2-1}{12}+(\Delta^\no)''(1) \ \ \ \ (\mbox{for $a\geq \deg\Delta$}).$$
This reproves the vanishing of ${\mathcal D}_a$ in \ref{g9}(${\mathcal D}$) for
such polynomials (cf. also with  \ref{g8}($**$)).

(4) Although the Alexander polynomial $\Delta(f)$  of the algebraic
knot $(S^3,K_f)$ ($f$ irreducible plane curve singularity) is known since
1932 \cite{Bu,Za}, and it was studied intensively (see e.g. \cite{Le,AC,EN}),
the property \ref{a2} remained hidden (to the best of the author's
knowledge).

On the other hand, similar properties were intensively studied in number
theory: namely in 40's, 50's and 60's a considerable
large  number of articles were published
about the coefficients of cyclotomic polynomials. Here we mention only a few
results.  If $\phi_n$ denotes the
$n^{th}$-cyclotomic polynomial, then it was proved that the coefficients of
$\phi_n$ have values in $\{-1,0,+1\}$ for $n=2^\alpha p^{\beta}q^{\gamma}$
($p$ and $q $ distinct odd primes) (result which goes back to the work of I.
Schur); if $n$ is a product of three distinct primes $pqr$ ($p<q<r$ and $
p+q>r$) then the coefficient of $t^r$ in $\phi_n$ is $-2$ (result of V.
Ivanov);
later Erd\H{o}s proved interesting estimates for the growth of the
coefficients; and G.S.
Kazandzidis provided exact formula for them. The interested
reader can consult \cite{NT}, pages 404-411, for a large list of articles
about this subject. (Reading these reviews, apparently the alternating property
was not perceived in this area either.)

Clearly, the above facts are not independent of our problem:
by \ref{a1}(3) the Alexander polynomial $\Delta(f)$  is a {\em product} of
cyclotomic polynomials.

(5) In fact, there is a recent result \cite{GDC}
in the theory of singularities,
which  implies the alternating property \ref{a2}(b). For any
irreducible curve singularity, using its normalization,
one can define a semigroup $S\subset\N$, with $0\in S$ and $\N\setminus S$
finite. Then,  in \cite{GDC}, based on some results of Zariski,
for an irreducible  plane curve singularities   is proved that
$\Delta(f)(t)/(1-t)=\sum_{i\in S}t^i$. This clearly implies the alternating
property.

(6) Are the irreducible plane curve singularities unique with the alternating
property? The answer is negative. In order to see this, consider the
Seifert integral homology sphere $\Sigma=\Sigma(a_1,a_2,\ldots, a_{k+1})$
(where  $\{a_i\}_i$ are pairwise coprime integers). Let $K$ be the special
orbit associated with the last arm (with Seifert invariant $a_{k+1}$).
Then the Alexander polynomial $\Delta_{\Sigma}(K)$ has alternating
coefficients. Indeed,
write $a:=a_1\cdots a_k$,  $a_i':=a/a_i$ for any $1\leq i\leq k$, and
let $S\subset \N$ be the semigroup (with $0\in S$)
generated by $a_1',\cdots, a_k'$. Then
$$\Delta_{\Sigma}(K)(t)=\frac{(1-t^a)^{k-1}(1-t)}
{\prod_{i=1}^k (1-t^{a_i'})}=
(1-t) \sum_{i\in S}t^i.$$
The first equality follows e.g. by \cite{EN}, section 11; the second by an
induction over $k$. This implies the alternating property as above.

(7) We can ask the following natural  question: what is that
property which distinguishes $(M,K_f)$
(where $f$ is an irreducible  plane curve singularity), or $(\Sigma,K)$
given  in (6),   from
the example described in  \ref{dnotzero} ? Why is in the first case the
${\mathcal D}$-invariant zero and in the second case not ? Can this be
connected with some property of the
semigroups $S$ associated with the curve whose link is $K$ ?
(We believe that the validity of an Abhyankar-Azevedo type theorem for this
curve  plays an important role in this phenomenon.)

(8) Examples show that the assumptions of \ref{a3} are really essential
(cf. also with \ref{rem}(2)).

\section{The link of $\{f(x,y)+z^n=0\}$}

\subsection{Preliminaries.}\label{l1} \
The present section is more technical than the previous ones, and
some of the details are left to the reader, which might cost the reader
some work.

Fix an irreducible plane curve
singularity $f:(\C^2,0)\to (\C,0)$ and let $K_f\subset S^3$ be its link
as in the previous section. Fix an integer $n\geq 1$, and consider the
``suspension'' germ $g:(\C^3,0)\to (\C,0)$ given by $g(x,y,z)=f(x,y)+z^n$.
Its link (i.e. $\{g=0\}\cap S^5_{\epsilon}$ for $\epsilon\ll 1$)
will be denoted by $M$. {\em
We will assume that $M$ is a rational homology sphere},
cf. \ref{pr}(c).

First, we recall/fix some numerical notations. We set:

\vspace{2mm}

\noindent
$\bullet $ \ the Newton pairs $\{(p_k,q_k)\}_{k=1}^s$ of $f$;\\
$\bullet$  \ the integers $\{a_k\}_{k=1}^s$ defined as in \ref{a1}(1);
recall that gcd$(p_k,a_k)=1$ for any $k$; \\
$\bullet$ \  $d_k:=\mbox{gcd}(n,p_{k+1}p_{k+2}\cdots p_s)$ for
$0\leq k\leq s-1$, and $d_s:=1$; \\
$\bullet$ \ $h_k:=d_{k-1}/d_k=\mbox{gcd}(p_k, n/d_k)$ for $1\leq k\leq s$;\\
$\bullet$ \ $\tilde{h}_k:=\mbox{gcd}(a_k,n/d_k)$ for $1\leq k\leq s$.

\vspace{2mm}

For any integer $1\leq l\leq s$,  let $M_{(l)}$ be the link of the
suspension singularity $g_{(l)}(x,y,z):=f_{(l)}(x,y)+z^{n/d_l}$. Evidently,
$M_{(s)}=M$, and $M_{(1)}=\Sigma(p_1,a_1,n/d_1)$.

\subsection{Some properties of the 3-manifolds  $\{M_{(l)}\}_l$.}\label{pr} \

{\bf (a)} \ \cite{DK,K,Ncy,Ste} \ For each $1\leq l\leq s$,
$M_{(l)}$ is the $(n/d_l)$-cyclic cover of $S^3$ branched along
$\{f_{(l)}=0\}$. Let $K_{(l)}\subset M_{(l)}$
be the preimage of $\{f_{(l)}=0\}$ with respect to  this cover.

{\bf (b)} \  $(M_{(l)},K_{(l)})$ can be represented by
(a ``canonical'')
plumbing graph (or resolution graph) which is compatible with the above
cover. This is made explicitly in \cite{nemsignat} (based on the idea of
\cite{Lbook}), see also \ref{et} here.  Using this
one obtains the following inductive picture.

For any $2\leq l\leq s$, $M_{(l)}$ can be obtained by  splicing, as
it is described in section 4,  the 3-manifold $M_2=\Sigma(p_l,a_l,n/d_l)$
along $K_2=\{y=0\}$ with $h_l$ copies of $M_{(l-1)}$ along the link
$K_{(l-1)}$ (with the same splicing data $\{A^{(i)}\}_i$ as in section 4).
(In order to prove this,
one needs to determine the invariant $M_w$ used in  \cite{nemsignat};
this is done in \cite{NDedI}, in  the proof of (3.2).)

{\bf (c)} \ Part (b) assures that $M$ is a rational homology sphere if and only if
for each $1\leq l\leq s$ the Seifert 3-manifold $\Sigma(p_l,a_l,n/d_l)$ is a
rational homology sphere. Since gcd$(p_l,a_l)=1$, this is happening if and
only if
\[
(h_l-1)(\tilde{h}_l-1)=0,\;\;\mbox{ for any $l$ (cf. \cite{Bries} or \cite{nemsignat})}
\]

{\bf (d)} \ Using (b) and \ref{s5}, one has:
\begin{equation*}
|H_1(M_{(l)})|=|H_1(\Sigma(p_l,a_l,n/d_l))|+h_l\cdot |H_1(M_{(l-1)})|
\ \ \mbox{for any $2\leq l\leq s$}.
\end{equation*}
or
\begin{equation*}
|H_1(M)|=\sum_{l=1}^s\, d_l\cdot |H_1(\Sigma(p_l,a_l,n/d_l))|.
\end{equation*}
In fact, one can give a complete description of the group $H_1(M)$
and the character group $H_1(M)\hat{}$ using \ref{s5} (we will come back
 to this  in \ref{et}).

{\bf (e)} \ As a parallelism, let us recall some  similar  formulae for other
numerical invariants: let $\mu_{(l)}$, respectively $\sigma_{(l)}$,
be the Milnor number, respectively the signature of the Milnor fiber
of $g_{(l)}$. Similarly, let $\mu(p_l,a_l,n/d_l)$ and
$\sigma(p_l,a_l,n/d_l)$ be the Milnor number and the signature
of the Brieskorn singularity $x^{p_l}+y^{a_l}+z^{n/d_l}$.
Then, by \cite{NDedI}, for any $2\leq l\leq s$:
\begin{equation*}
\sigma_{(l)}=\sigma(p_l,a_l,n/d_l)+h_l\cdot \sigma_{(l-1)}\ \ \mbox{or} \ \
\sigma_{(s)}=\sum_{l=1}^s\, d_l\cdot \sigma(p_l,a_l,n/d_l).
\end{equation*}
By contrast, for the Milnor numbers one has
$\mu_{(l)}=\mu(p_l,a_l,n/d_l)+p_l\cdot \mu_{(l-1)}$
(involving $p_l$ versus $h_l$, fact
which follows e.g. from  \ref{a1}(4)).

Our goal is to establish an inductive formula for
$\ssw^0_M(\sigma_{can})$, similar to $|H_1|$ or to the signature
$\sigma_{(l)}$. (For $\lambda_W$ or $\et_{M,\sigma_{can}}(1)$ such  a
formula does not hold, see below.)

{\bf (f)} \ \cite{DK,K,Ncy,Ste} \   $(S^3,K_{f_{(l)}})$ is fibrable.
Let $\M_{geom,(l)}:F_{(l)}\to F_{(l)}$ (respectively $\M_{(l)}$)
be a geometric  (respectively the algebraic) monodromy acting on the
Milnor fiber $F_{(l)}$ (respectively on $H_1(F_{(l)}$).
Then $(M_{(l)},K_{(l)})$ is also fibrable,
whose open book decomposition has the same
fiber $F_{(l)}$ and  geometric monodromy $\M_{geom,(l)}^{n/d_l}$.
In particular, the (normalized) Alexander polynomial of $(M_{(l)},K_{(l)})$
is (the normalization of)
\[
\Delta_{M_{(l)}}(t)=\Delta_{M_{(l)}}(K_{(l)})(t)= \det(1-t\M^{n/d_l}_{(l)}).
\]
Therefore, using \ref{s0e}(7)  and $2r_{l}:=\mbox{rank}H_1(F_{(l)})$:
$$\Delta^{\no}_{M_{(l)}}(t)=\frac{1}{|H_1(M_{(l)})|}\cdot t^{-r_{l}}\cdot
\det(1-t\M^{n/d_l}_{(l)})
\ \ \mbox{with} \ \
|H_1(M_{(l)})|=|\det(1-\M^{n/d_l}_{(l)})|.
$$
Notice that $\Delta^{\no}_{M_{(l)}}(t)$ can be deduced from the Alexander
polynomial
\[
\Delta(f_{(l)})(t)=\det(1-t\M_{(l)}).
\]
of $f_{(l)}$ (cf. section 5). Indeed, for any
polynomial $\Delta(t)$ of degree $2r$ and of the form
$\Delta(t)=\prod_v\, (1-t^{m_v})^{n_v}$,
and for any positive  integer $k$, define
\[
\Delta^{c(k)}(t):= \prod_v\,   (1-t^{m_v/gcd(m_v,k)})^{n_v\cdot
gcd(m_v,k)}.
\]
Let $\Delta^{c(k),\no}(t)=t^{-r}\Delta^{c(k)}(t)/\Delta^{c(k)}(1)$ denote
the normalization of $\Delta^{c(k)}(t)$. An
eigenvalue-argument then proves:
$$\Delta_{M_{(l)}}(t)=\Delta(f_{(l)})^{c(n/d_l)}(t)\ \ \mbox{and}\ \
\Delta^{\no}_{M_{(l)}}(t)=\Delta(f_{(l)})^{c(n/d_l),\no}(t).$$

{\bf (g)} \ The inductive formula \ref{a1}(4) reduces the computation
of the Alexander invariants to the Seifert case.   Clearly
(from \ref{a1}(5), \ref{g2}(c) and (f) above):
$$\Delta(x^{p_l}+y^{a_l})^{c(n/d_l),\no}(t)=\frac{1}{a_l^{h_l-1}
p_l^{\tilde{h}-1}}\cdot
t^{-(a_l-1)(p_l-1)/2}\cdot \frac{(t^{p_la_l/(h_l\tilde{h}_l)}-1)
^{h_l\tilde{h}_l}(t-1)}{(t^{p_l/h_l}-1)^{h_l}
(t^{a_l/\tilde{h}_l}-1)^{\tilde{h}_l}}.$$
(Recall that $(h_l-1)(\tilde{h}_l-1)=0$ for any $l$.)
Then, by a computation,  one can show that
$$(\Delta(x^{p_l}+y^{a_l})^{c(n/d_l),\no})''(1)=
\frac{1}{12} \Big(\frac{a_l^2}{\tilde{h}_l}-1\Big)\Big(\frac{p_l^2}{h_l}
-1\Big).$$

{\bf (h)} \ Using  (f) and \ref{a1}(4) one gets:
$$\Delta_{M_{(l)}}^{\no}(t)=
\Delta(x^{p_l}+y^{a_l})^{c(n/d_l),\no}(t)\cdot
\Big[\Delta_{M_{(l-1)}}^{\no}(t^{p_l/h_l})\Big]^{h_l}.$$
Then, using $(\Delta^{\no})'(1)=0$ (cf. \ref{a4}(2)) and the result from (g),
 one obtains:
\begin{equation*}
 (\Delta_{M_{(l)}}^\no)''(1)=
\frac{1}{12} \Big(\frac{a_l^2}{\tilde{h}_l}-1\Big)\Big(\frac{p_l^2}{h_l}
-1\Big)
+\frac{p_l^2}{h_l}\cdot (\Delta_{M_{(l-1)}}^\no)''(1).
\end{equation*}
Therefore:
\begin{equation*}
(\Delta_{M_{(l)}}^\no)''(1)=\sum_{k=1}^l\, \frac{1}{12}\,
 \Big(\frac{a_k^2}{\tilde{h}_k}-1\Big)\Big(\frac{p_k^2}{h_k}-1\Big)\cdot
\frac{(p_{k+1}\cdots p_l)^2}{h_{k+1}\cdots h_l}.
\end{equation*}

{\bf (i)} \ For any $2\leq l\leq s$, one has:
$$\lambda_W(M_{(l)})=\lambda_W(\Sigma(p_l,a_l,n/d_l))+h_l\cdot
\lambda_W(M_{(l-1)})\hspace{2cm}$$ $$\hspace{2.8cm}
+\frac{np_l(h_l-1)}{d_la_lh_l}\cdot
\sum_{k=1}^{l-1}\, \frac{1}{12}\,
 \Big(\frac{a_k^2}{\tilde{h}_k}-1\Big)\Big(\frac{p_k^2}{h_k}-1\Big)\cdot
\frac{(p_{k+1}\cdots p_{l-1})^2}{h_{k+1}\cdots h_{l-1}}.$$
Indeed, if $h_l=1$ then we have additivity as in \ref{s2}(4), cf. also
with the proof of \ref{swc}, part B. If $\tilde{h}_l=1$, then apply
\ref{g7} (see also \ref{g9}) and (h) above.

For the value of the
Casson-Walker invariant $\lambda_W(\Sigma(p,a,n))$ of a Seifert
manifold, see \cite[6.1.1]{Lescop}  or \cite[5.4]{NN}.

{\bf (j)} \ Below in (k),  we will compute $(\Delta^\no_{M_{(l)}})''(1)$
in terms of $\{\, (\Delta(f_{(k)})^\no)''(1)\, \}_{k\leq l}$.
Clearly, one can obtain similar  inductive formula for these
$(\Delta(f_{(k)})^\no)''(1)$ as that one in (h) by taking $n=1$.
More precisely, for any $2\leq l\leq s$:
\begin{equation*}
 (\Delta(f_{(l)})^\no)''(1)=
(a_l^2-1)(p_l^2-1)/12 + p_l^2\cdot (\Delta(f_{(l-1)})^\no)''(1).
\end{equation*}

{\bf (k)} \ The next (rather
complicated) identity looks very artificial, but it is
one of the most important formulae  in this list. Basically, it validity
is equivalent
with the fact that the two correction terms ${\mathcal O}(\et)$ and
${\mathcal O}(\lambda_W/2)$ are the same, cf. \ref{splet} and \ref{swadd}.
 For any $2\leq l\leq s$ one has:
\begin{equation*}
(\Delta^\no_{M_{(l)}})''(1)=
(\Delta(f_{(l)})^\no)''(1)
-\sum _{k=1}^l\, \frac{a_k^2p_k^2}{\tilde{h}_k^2h_k^2}
\cdot\frac{p_{k+1}^2\cdots p_l^2}{h_{k+1}\cdots h_l}\cdot A_k,
\end{equation*}
where
\begin{equation*}
A_k:= \frac{h_k(h_k-1)}{a_k^2}\cdot
\Big[ \,  \frac{a_k^2-1}{12} +(\Delta(f_{(k-1)})^\no)''(1)\, \Big]+
 \frac{\tilde{h}_k(\tilde{h}_k-1)}{p_k^2}\cdot \frac{p_k^2-1}{12}.
\end{equation*}
For the proof proceed as follows. Let $E_{(l)}$ be the difference between
the left  and the right hand side of the identity. Then, using
the inductive formulae (h) and (j) and the property $(h_k-1)(\tilde{h}_k-1)=0$,
 by an elementary computation one can verify that
$E_{(l)}-(p_l^2/h_l)\cdot E_{(l-1)}=0$, and $E_{(1)}=0$. Then $ E_{(l)}=0$ by
induction.

{\bf (l)} The last invariant we wish to determine
is $\et_{M_{(l)},\sigma_{can}}(1)$. The
computation is more involved and it is separated in the next subsections.
The inductive formula for $\et_{M_{(l)},\sigma_{can}}(1)$ is
given in \ref{splet}.

\subsection{Characters of $H_1(M)$.}\label{et} \
In the computation of  $\et_{M,\sigma_{can}}(1)$
we plan to use \ref{sw}(16). For this, we need to describe the characters
$\chi\in \hat{H}=H_1(M)\hat{}$.

The group $H$ can be  determined in many different ways. For example,
using the monodromy operator $\M=\M_{(s)}$ of $f$,
$H$ can be identified, as an abstract
group,  with coker$(1-\M^n)$. The homology of the Milnor fiber of $f$ and
$\M$ have  a direct sum decomposition
with respect to the splicing (see e.g. \cite{EN} or \cite{NDedI}).
This can be used to provide an inductive description of $H$.

Nevertheless, we prefer to use \ref{s5}.
The main reason is that, in fact, we have to understand  $\hat{H}$ (rather
than $H$) together  with the description of the supports
$\{v : \ \chi(g_v)\not=1\}$ for each character $\chi\in \hat{H}$. For this,
the discussion from section 3 is more suitable.

We consider again the ``canonical'' plumbing graph $\Gamma(M)$
of  $M$ provided by the algorithm \cite{nemsignat}, cf. \ref{pr}(b) here.
In fact, that algorithm
provides $\Gamma:=\Gamma(M,K_z)$, the plumbing graph of the 3-dimensional
link $M=\{f(x,y)+z^n=0\}$
 with the knot  $K_z:=\{z=0\}$ in it. Again, if we replace  the strings by dash-lines, then one can represent
$\Gamma$ as a covering graph of $\Gamma(S^3,K_f)$; for details, see
[loc. cit.], cf.  also with  \ref{pr}(b).
If we denote this graph-projection by $\pi$, then
\[
\#\pi^{-1}(v_k)=h_{k+1}\cdots h_s, \;\;1\leq k\leq s,
\]
\[
\ \ \ \ \
\#\pi^{-1}(\bar{v}_k)=\tilde{h}_kh_{k+1}\cdots h_s, \;\; 1\leq k\leq s,
\ \ \mbox{and}\]
\[
\#\pi^{-1}(\bar{v}_0)=h_1\cdots h_s
\]
(see \ref{a1} for notations about $\Gamma(S^3,K_f)$).
In fact, there is
a $\Z_n$-action on $\Gamma$ which acts transitively on each fiber of $\pi$,
hence all the vertices above a given vertex $v\in {\mathcal V}^*$ of
$\Gamma(S^3,K_f)$  are symmetric in $\Gamma$.  In particular, their
decorations and their numerical invariants (computed out of the graph $\Gamma$)
are the same.
Therefore, $\Gamma$ has the following schematic form:

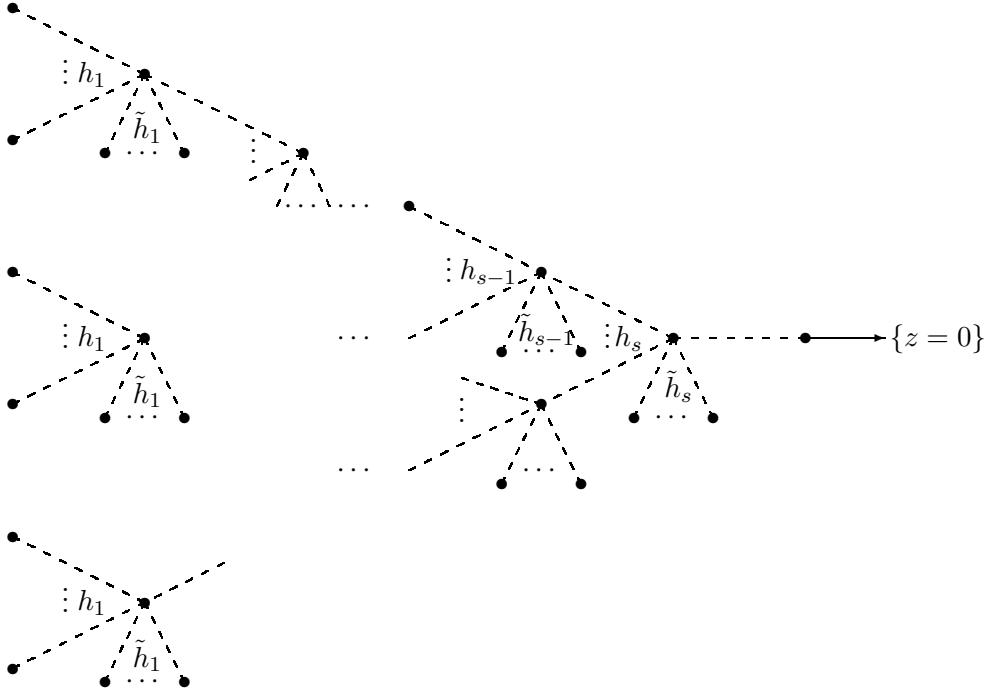
\begin{picture}(400,300)(30,-50)
\put(50,225){\circle*{4}}
\put(50,175){\circle*{4}}
\put(50,125){\circle*{4}}
\put(50,75){\circle*{4}}
\put(50,25){\circle*{4}}
\put(50,-25){\circle*{4}}
\put(100,200){\circle*{4}}
\put(100,100){\circle*{4}}
\put(100,0){\circle*{4}}
\dashline{3}(100,200)(50,225)
\dashline{3}(100,200)(50,175)
\dashline{3}(100,100)(50,125)
\dashline{3}(100,100)(50,75)
\dashline{3}(100,0)(50,25)
\dashline{3}(130,15)(50,-25)
\put(85,170){\circle*{4}}
\put(85,70){\circle*{4}}
\put(85,-30){\circle*{4}}
\put(115,170){\circle*{4}}
\put(115,70){\circle*{4}}
\put(115,-30){\circle*{4}}
\dashline{3}(100,200)(85,170)
\dashline{3}(100,100)(85,70)
\dashline{3}(100,0)(85,-30)
\dashline{3}(100,200)(115,170)
\dashline{3}(100,100)(115,70)
\dashline{3}(100,0)(115,-30)
\put(160,170){\circle*{4}}
\dashline{3}(160,170)(100,200)
\dashline{3}(160,170)(140,160)
\dashline{3}(160,170)(150,150)
\dashline{3}(160,170)(170,150)
\put(200,150){\circle*{4}}
\put(250,125){\circle*{4}}
\put(300,100){\circle*{4}}
\put(250,75){\circle*{4}}
\put(350,100){\circle*{4}}
\dashline{3}(200,150)(300,100)
\dashline{3}(200,50)(300,100)
\dashline{3}(250,125)(200,100)
\dashline{3}(250,125)(235,95)
\dashline{3}(250,125)(265,95)
\dashline{3}(250,75)(235,45)
\dashline{3}(250,75)(265,45)
\put(235,45){\circle*{4}}
\put(235,95){\circle*{4}}
\put(265,95){\circle*{4}}
\put(265,45){\circle*{4}}
\dashline{3}(300,100)(285,70)
\dashline{3}(300,100)(315,70)
\put(285,70){\circle*{4}}
\put(315,70){\circle*{4}}
\dashline{3}(300,100)(350,100)
\dashline{3}(250,75)(220,85)
\put(400,100){\makebox(0,0){$\{z=0\}$}}
\put(350,100){\vector(1,0){30}}
\put(70,204){\makebox(0,0){$\vdots$}}
\put(80,200){\makebox(0,0){$h_1$}}
\put(70,104){\makebox(0,0){$\vdots$}}
\put(80,100){\makebox(0,0){$h_1$}}
\put(70,4){\makebox(0,0){$\vdots$}}
\put(80,0){\makebox(0,0){$h_1$}}
\put(100,170){\makebox(0,0){$\cdots$}}
\put(101,180){\makebox(0,0){$\tilde{h}_1$}}
\put(100,70){\makebox(0,0){$\cdots$}}
\put(101,80){\makebox(0,0){$\tilde{h}_1$}}
\put(100,-30){\makebox(0,0){$\cdots$}}
\put(101,-20){\makebox(0,0){$\tilde{h}_1$}}
\put(141,174){\makebox(0,0){$\vdots$}}
\put(215,129){\makebox(0,0){$\vdots$}}
\put(230,125){\makebox(0,0){$h_{s-1}$}}
\put(220,76){\makebox(0,0){$\vdots$}}
\put(275,104){\makebox(0,0){$\vdots$}}
\put(283,100){\makebox(0,0){$h_{s}$}}
\put(250,50){\makebox(0,0){$\cdots$}}
\put(250,95){\makebox(0,0){$\cdots$}}
\put(252,103){\makebox(0,0){$\tilde{h}_{s-1}$}}
\put(300,70){\makebox(0,0){$\cdots$}}
\put(302,82){\makebox(0,0){$\tilde{h}_{s}$}}
\put(180,150){\makebox(0,0){$\cdots$}}
\put(160,150){\makebox(0,0){$\cdots$}}
\put(180,100){\makebox(0,0){$\cdots$}}
\put(180,50){\makebox(0,0){$\cdots$}}
\end{picture}

Now, we consider the Seifert manifold $\Sigma:=\Sigma(p,a,n)$ with
gcd$(p,a)=1$.

If gcd$(a,n)=1$ and gcd$(p,n)=d$, by  \ref{g3}, $H_1(\Sigma)$
is generated by the homology classes of $\{K_2^{(i)}\}_{i=1}^d$. For an
arbitrary  character $\chi$, we write $\chi(K_2^{(i)})=\xi_i$. Here
$\xi_i$ is an $a$-root of unity in $\C$ (shortly $\xi_i\in \Z_a$).
Then $\chi$ is completely characterized by the collection $\{\xi_i\}
_{i=1}^d$ which satisfy $\xi_i\in \Z_a$ for any $i$ and $\xi_1\cdots \xi_d=1$.

Notice that $\chi(O)=1$, and $\chi$ is supported by those $d$ ``arms'' of the
 star-shaped graph which have Seifert invariant $a$ (i.e. $\chi(g_v)=1$
for the vertices $v$ situated  on the other arms).
In \cite{NN} is proved that for each
non-trivial character $\chi$, in $\hat{P}_{\Sigma,\chi,u}(t)$ one can take
for $u$ the central vertex $O$. Moreover, $\lim_{t\to 1}
\hat{P}_{\Sigma,\chi,u}(t)$ can be non-zero only if $\chi$ is supported
exactly on two arms,  i.e. $\xi_i\not=1$ exactly for two
values of $i$, say for $i_1$ and $i_2$ (hence $\xi_{i_1}=\bar{\xi}_{i_2}$).

If gcd$(n,p)=1$,  but gcd$(n,a)\not=1$, then clearly we have a symmetric
situation; in this case  we use the notation $\eta$ instead of $\xi$.

In both situation, for any character $\chi$, $\chi(g_v)=1$ for any
$v$ situated on the arm with Seifert invariant $n/\mbox{gcd}(n,ap)$.

These properties proved for the building block $\Sigma$ will generate all
the properties of $H=H_1(M)$ via the splicing properties \ref{s5}
and linking relations \ref{g3} (L\#).

For example,  one can prove by induction, that for any character
$\chi$ of $H$, $\chi(g_v)=1$ for any vertex $v$ situated on the string
which supports the arrow of $K_z$.

Consider the splicing decomposition
\[
M=M_{(s)}=h_sM_{(s-1)}\coprod \Sigma(p_s,a_s,n).
\]
 As in the previous inductive arguments, assume that we understand the characters $\chi$ of $H_1(M_{(s-1)})$.
By \ref{s5}, they can be considered in a natural way as characters of
$H_1(M_{(s)})$ satisfying additionally  $\chi(g_v)=1$
for any  vertex $v$ of $\Sigma(p_s,a_s,n)$ (cf. also with the first paragraph
of the proof of \ref{s7}).
We say that these characters do not ``propagate'' from
$M_{(s-1)}$ into $\Sigma(p_s,a_s,n)$.

In the ``easy'' case when $h_s=1$ (even if $\tilde{h}_s\not=1$),
the splicing invariants are
$o_1=o_2=1$ and $k_1=k_2=0$, hence $H$ (together with its linking form)
is a direct sum in a natural way,
hence the characters of $\Sigma(p_s,a_s,n)$ (described above)
will not propagate into $M_{(s-1)}$ either.

On the other hand, if $h_s>1$, then the {\em non-trivial characters}
 of $\Sigma(p_s,a_s,n)$
do propagate into the $h_s$ copies of $M_{(s-1)}$.
If $\chi$ is such a character, with $\chi(K_2^{(i)})=\xi_i$, then
analyzing the properties of the linking numbers in \ref{s5} we deduce that
$\chi$ can be extended into  $M$ in such a way
 that $\chi(g_{v^{(i)}})=\xi_i^l$  for any
vertex $v^{(i)}$ of  $M^{(i)}_{(s-1)}$, where
$l=Lk_{M^{(i)}_{(s-1)}}(g_{v^{(i)}},K_1^{(i)})$.
This fact can be proved as follows.
Assume that $\chi=\exp(Lk_\Sigma(L,\cdot))$ for some fixed $L\subset \Sigma$.
Then $\xi_i=\chi(K_2^{(i)})=\exp(Lk_\Sigma(L,K_2^{(i)}))$. Therefore,
$\chi$ extended into $M$ as $\exp(Lk_M(L,\cdot))$
satisfies (cf. \ref{s5}(11))
\[
\chi(g_{v^{(i)}})=\exp(Lk_M(L,g_{v^{(i)}}))=\exp(Lk_\Sigma(L,K_2^{(i)})\cdot Lk_{M_{s-1}^{(i)}}(K_1^{(i)},g_{v^{(i)}}))=\xi_i^l.
\]
Clearly, all these linking numbers
$Lk_{M^{(i)}_{(s-1)}}(g_{v^{(i)}},K_1^{(i)})$
 can be determined inductively by the formulae provided in
\ref{s5} and \ref{g3} (L\#).

Moreover, if we multiply  such an extended  character $\chi$  with a character
$\chi'$ supported by $M_{(s-1)}$, then $\chi\chi'(g_v)=\chi(g_v)$ for
any vertex $v$ situated above the vertices $v_s$ or $v_{s-1}$ or above
any vertex  on the edge connecting $v_{s-1}$ with $v_s$
(since the support of $\chi'$ does not contain these vertices).
 We will say that such a  character $\chi\chi'$  is ``born at level $s$''
(provided that the original $\chi$ of $\Sigma$  was  non-trivial).

(The interested reader can reformulate the above discussion in the
language of an exact sequence of dual groups, similar to one used
in the first paragraph of the proof of \ref{s7}.)

\vspace{2mm}

Below we provide some examples. In this diagrams,
for any fixed character $\chi$, we put on the vertex
$v$ the complex number $\chi(g_v)$.

\subsection{Example.}\label{seq2} Assume that $s=2$. Then basically one has two
different cases: $\tilde{h}_1=h_2=1$, or $\tilde{h}_1=\tilde{h}_2=1$ (since
the first pairs $(p_1,a_1)$ can be permuted).

In the first (``easy'') case, the schematic diagram of the  characters is:

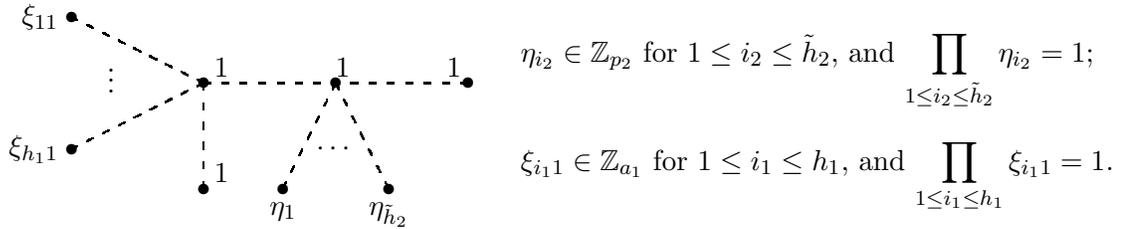
\begin{picture}(400,110)(170,60)
\put(200,150){\circle*{4}}
\put(200,100){\circle*{4}}
\put(250,125){\circle*{4}}
\put(300,125){\circle*{4}}
\put(350,125){\circle*{4}}
\dashline{3}(200,150)(250,125)
\dashline{3}(250,125)(200,100)
\dashline{3}(250,125)(350,125)
\dashline{3}(250,125)(250,85)
\put(250,85){\circle*{4}}
\put(280,85){\circle*{4}}
\put(320,85){\circle*{4}}
\dashline{3}(300,125)(280,85)
\dashline{3}(300,125)(320,85)
\put(215,129){\makebox(0,0){$\vdots$}}
\put(300,100){\makebox(0,0){$\cdots$}}
\put(254,131){\makebox(0,0)[l]{$1$}}
\put(254,91){\makebox(0,0)[l]{$1$}}
\put(300,131){\makebox(0,0)[l]{$1$}}
\put(342,131){\makebox(0,0)[l]{$1$}}
\put(280,80){\makebox(0,0)[t]{$\eta_1$}}
\put(320,80){\makebox(0,0)[t]{$\eta_{\tilde{h}_2}$}}
\put(195,150){\makebox(0,0)[r]{$\xi_{11}$}}
\put(195,100){\makebox(0,0)[r]{$\xi_{h_11}$}}
\put(370,130){\makebox(0,0)[l]{$\eta_{i_2}\in \Z_{p_2}$ for $1\leq i_2\leq
\tilde{h}_2$, and $ \displaystyle{\prod_{1\leq i_2\leq \tilde{h}_2}\eta_{i_2}=1}$;}}
\put(370,90){\makebox(0,0)[l]{$\xi_{i_11}\in \Z_{a_1}$ for $1\leq i_1\leq h_1$, and $ \displaystyle{\prod_{1\leq i_1\leq h_1}\xi_{i_11}=1}$.}}
\end{picture}

In the second case $\tilde{h}_1=\tilde{h}_2=1$, with the notation
$p_1':=p_1/h_1$, the characters $\hat{H}$ are as follows.

\begin{picture}(400,150)(170,30)
\put(200,150){\circle*{4}}
\put(250,125){\circle*{4}}
\put(300,100){\circle*{4}}
\put(250,75){\circle*{4}}
\put(350,100){\circle*{4}}
\dashline{3}(200,150)(300,100)
\dashline{3}(200,50)(300,100)
\dashline{3}(250,125)(200,110)
\put(200,110){\circle*{4}}
\put(200,90){\circle*{4}}
\put(200,50){\circle*{4}}
\dashline{3}(250,125)(250,95)
\dashline{3}(250,75)(250,45)
\put(250,45){\circle*{4}}
\put(250,95){\circle*{4}}
\dashline{3}(300,100)(300,70)
\put(300,70){\circle*{4}}
\dashline{3}(300,100)(350,100)
\dashline{3}(250,75)(200,90)
\put(215,129){\makebox(0,0){$\vdots$}}
\put(220,76){\makebox(0,0){$\vdots$}}
\put(275,104){\makebox(0,0){$\vdots$}}

\put(258,129){\makebox(0,0)[l]{$\xi_1^{a_1p_1'}$}}
\put(258,98){\makebox(0,0)[l]{$\xi_1^{a_1}$}}
\put(258,70){\makebox(0,0)[l]{$\xi_{h_2}^{a_1p_1'}$}}
\put(258,45){\makebox(0,0)[l]{$\xi_{h_2}^{a_1}$}}
\put(304,70){\makebox(0,0)[l]{$1$}}
\put(304,106){\makebox(0,0)[l]{$1$}}
\put(342,106){\makebox(0,0)[l]{$1$}}

\put(195,150){\makebox(0,0)[r]{$\xi_{11}\xi_1^{p_1'}$}}
\put(195,110){\makebox(0,0)[r]{$\xi_{h_11}\xi_1^{p_1'}$}}
\put(195,90){\makebox(0,0)[r]{$\xi_{1h_2}\xi_{h_2}^{p_1'}$}}
\put(195,50){\makebox(0,0)[r]{$\xi_{h_1h_2}\xi_{h_2}^{p_1'}$}}

\put(370,130){\makebox(0,0)[l]{$\xi_{i_2}\in \Z_{a_2}$ for $1\leq i_2\leq h_2$, and $ \displaystyle{\prod_{1\leq i_2\leq h_2}\xi_{i_2}=1}$;}}
\put(400,100){\makebox(0,0)[l]{and for any fixed $i_2$:}}
\put(370,70){\makebox(0,0)[l]{$\xi_{i_1i_2}\in \Z_{a_1}$ for $1\leq i_1\leq h_1$, and $ \displaystyle{\prod_{1\leq i_1\leq h_1}\xi_{i_1i_2}=1}$.}}
\end{picture}

\subsection{\ }\label{ww} Using the discussion \ref{et},
the above example for  $s=2$ can be generalized inductively
to arbitrary $s$. For the convenience of the reader, we make this explicit.
In order to have a uniform notation, we consider the case
$\tilde{h}_k=1$ for any $1\leq k\leq s$.
The interested reader is invited to write down similar description of
the characters in those cases when $\tilde{h}_k\not=1$ for some $k$,
using the present model and \ref{seq2}.

More precisely, for any character $\chi$,  we will
 indicate $\chi(g_{v'})$ for any vertex $v'\in \pi^{-1} ({\mathcal V}^*)$.

It is convenient to introduce the index set $(i_1,i_2,\ldots, i_s)$,
where $1\leq i_l\leq h_l$ for any $1\leq l\leq s$.
As we already mentioned, this set can be considered as the index set of
$\pi^{-1}(\bar{v}_0)$. Moreover, for any $1\leq k\leq s-1$,
$(i_{k+1},\ldots,i_s)$ (where  $1\leq i_l\leq h_l$ for any $k+1\leq l\leq s$)
is the index set of $\pi^{-1}(v_k)$ (respectively of $\pi^{-1}(\bar{v}_k)$
since $\tilde{h}_k=1$). Moreover, for any $l$ we write $p_l':=p_l/h_l$.

Next, we consider a system of roots of unity as follows:

\vspace{2mm}

\noindent
$\bullet$ \ $\xi_{i_s}\in \Z_{a_s}$ with
$\prod_{1\leq i_s\leq h_s}\, \xi_{i_s}=1$;\\
$\bullet$ \ for any fixed $i_s$ a collection $\xi_{i_{s-1}i_s}\in \Z_{a_{s-1}}$
with
$\prod_{1\leq i_{s-1}\leq h_{s-1}}\, \xi_{i_{s-1}i_s}=1$;\\
and, more generally, if $1\leq k\leq s-1$:\\
$\bullet$ \ for any fixed $(i_{k+1},\,\cdots\, , i_s)$ a collection
$\xi_{i_k\,\cdots\, i_s}\in \Z_{a_{k}}$
with
$\prod_{1\leq i_{k}\leq h_{k}}\, \xi_{i_k\,\cdots\, i_s}=1$.

\vspace{2mm}

\noindent Then,  any character
$\chi$  can be characterized by the following properties:

\vspace{2mm}

\noindent $\bullet$ \
$\pi^{-1}(v_s)$ contains exactly one vertex, say $v'$. Then $\chi(g_{v'})=1$.
The same is valid for $\pi^{-1}(\bar{v}_s)$ and for that vertex in $\Gamma$
which supports the arrow $\{z=0\}$.\\
$\bullet$ \
For any $1\leq k\leq s-1$, if $v'(i_{k+1},\ldots,i_s)$ is the vertex in
$\pi^{-1}(v_k)$ corresponding to the index $(i_{k+1},\ldots,i_s)$, then
$$\chi(g_{v'(i_{k+1},\ldots,i_s)})=
\xi_{i_{k+1}\,\cdots\, i_s}^{a_k\,p_k'}\cdot \,\xi_{i_{k+2}\,\cdots\, i_s}
^{a_k\,p_k'\, p_{k+1}'}\, \cdots\, \xi_{i_s}^{a_k\,p_k'\,\cdots \,p_{s-1}'}.$$
$\bullet$ \ Similarly,
 if $\bar{v}'(i_{k+1},\ldots,i_s)$ is the vertex in
$\pi^{-1}(\bar{v}_k)$ corresponding to the index $(i_{k+1},\ldots,i_s)$, then
$$\chi(g_{\bar{v}'(i_{k+1},\ldots,i_s)})=
\xi_{i_{k+1}\,\cdots\, i_s}^{a_k}\,\cdot \xi_{i_{k+2}\,\cdots\, i_s}
^{a_k\, p_{k+1}'}\, \cdots\, \xi_{i_s}^{a_k\,p_{k+1}'\,\cdots \,p_{s-1}'}.$$
$\bullet$ \ Finally,  if
$\bar{v}'(i_{1},\ldots,i_s)$ is the vertex in
$\pi^{-1}(\bar{v}_0)$ corresponding to the index $(i_{1},\ldots,i_s)$, then
$$\chi(g_{\bar{v}'(i_{1},\ldots,i_s)})=
\xi_{i_{1}\,\cdots\, i_s} \cdot \xi_{i_{2}\,\cdots\, i_s}
^{p_1'}\, \cdots\, \xi_{i_s}^{p_1'\,\cdots \,p_{s-1}'}.$$

If $\xi_{i_s}\not=1$ for some $i_s$ then $\chi$ is ``born at level $s$''.
If $\xi_{i_s}=1$ for all $i_s$, but $\chi_{i_{s-1}i_s}\not=1$ for some
$(i_{s-1},i_s)$, then $\chi$ is ``born at level $s-1$'', etc.

In general, a character $\chi$ is ``born at level $k$'' ($1\leq k\leq s$)
if for any $l\geq k$ and $v'\in \pi^{-1}(v_l)$, one has $\chi(g_{v'})=1$,
but there exists at least one vertex $v'\in \pi^{-1}(v_k)$ which is adjacent
(in the graph $\Gamma$) with the support of $\chi$.

\vspace{2mm}

The next result  will be crucial when we apply the Fourier inversion
formula \ref{sw}(14). It is a really remarkable property of the links
associated with irreducible plane curve singularities. It is the most
important qualitative ingredient in our torsion computation (see also
\ref{rem} for another powerful application).

\subsection{Proposition.}\label{sum} {\em Consider $(S^3,K_f)$
associated with an irreducible plane curve singularity $f$. Let
$\Delta_{S^3}(K_f)(t)$ be its Alexander polynomial. For any integer $n\geq 1$,
consider $(M,K_z)$, i.e. the link $M$ of $\{f(x,y)+z^n=0\}$  and the knot
$K_z:=\{z=0\}$ in it. Let $\Delta^H_M(K_z)(t)$ be the Alexander invariant
defined in \ref{s0e} (10-11) with $H=H_1(M)$. Then }
$$\Delta^H_M(K_z)(t)=\Delta_{S^3}(K_f)(t).$$
\begin{proof}
First we notice that in \ref{s0e}(10), $g_u=K_z$ and $o(u)=1$. Then, by
\ref{s0d}(5), $m_{v'}(u)=Lk_M(g_{v'},K_z)$ for any vertex $v' $ of $\Gamma$.
By the algorithm \cite{nemsignat},
for any vertex $v'\in \pi^{-1}(v)$, where  $v\in
{\mathcal V}^*(\Gamma(S^3,K_f))$,  the linking number
$Lk_M(g_{v'},K_z)$ is given by $w_{v'}=w_v/\mbox{gcd}(w_v,n)$.
Recall, that for any $v\in {\mathcal V}^*(\Gamma(S^3,K_f))$, the corresponding
weights $w_v$ are given in \ref{a1}(2). In particular, this discussion provides
all the weight $w_{v'}(u)$ needed in the definition \ref{s0e}(10) of
$\Delta_{M,\chi}(K_z)(t)$. One the other hand,
for characters $\chi\in \hat{H}$ we can use the above description.
Then the proposition follows (by some computation) inductively using the
Algebraic Lemma, part (b), for the Alexander polynomials $\Delta(f_{(l)}(0))$
($0\leq l\leq s-1$, where $\Delta(f_{(0)})(t)\equiv 1$). Notice that this
lemma can be applied  thanks to the proposition \ref{a2} (which assures
that the coefficients of $\Delta(f_{(l)}(t)$ are alternating), and to the
inequality \ref{a1}(6) (which assures that $a$ ``is sufficiently large'').
 (For an expression of $\Delta(f_{(l)}(t)$, see \ref{a1}).)

In the next example we make this argument explicit
for the case $s=2$ and $\tilde{h}_1=\tilde{h}_2=1$. Using this model,
the reader can complete the general case easily.
\end{proof}

\subsection{Example.}\label{seq2b} Assume that  $s=2$ and
$\tilde{h}_1=\tilde{h}_2=1$. Then, with the notation of \ref{seq2},
 $\Delta_{M,\chi}(K_z)(t)$ equals
$$\prod_{i_1,i_2}\frac{1}{1-t^{p_1'\, p_2'}\xi_{i_1i_2}\, \xi_{i_2}^{p_1'}}
\cdot \prod_{i_2}
\frac{(1-t^{a_1\,p_1'\, p_2'}\xi_{i_2}^{a_1\, p_1'})^{h_1}}
{1-t^{a_1p_2'}\xi_{i_2}^{a_1}}\, \cdot\,
\frac{(1-t^{a_2\, p_2'})^{h_2}\cdot (1-t)}{1-t^{a_2}}.$$
First, for each fixed index $i_2$,  we make a sum over $\xi_{i_1i_2}\in
 \Z_{a_1}$. Using \ref{a3}(b) for $\Delta\equiv 1$, $t=t^{p_1'\, p_2'}\xi_{i_2}
^{p_1'}$ and $a=a_1$ and $d=h_1$, the above expression transforms (after some
simplifications)  into:
$$\prod_{i_2}
\frac{1-t^{a_1\,p_1\, p_2'}\xi_{i_2}^{a_1\, p_1}}
{(1-t^{p_1\, p_2'}\xi_{i_2}^{p_1})\cdot
(1-t^{a_1p_2'}\xi_{i_2}^{a_1})}\, \cdot\,
\frac{(1-t^{a_2\, p_2'})^{h_2}\cdot (1-t)}{1-t^{a_2}}.$$
The expression in the product is exactly $\Delta(f_{(1)})(t^{p_2'})/
(1-t^{p_2'})$. Therefore, \ref{a3}(b) can be applied again, now for
$\Delta=\Delta(f_{(1)})$, $t=t^{p_2'}$, $a=a_2$ and $d=h_2$. Then the
expression transforms into $\Delta(f_{(2)})$.

\subsection{Remarks.}\label{rem} \
(1) If the link $M$  of $\{f(x,y)+z^n=0\}$ is a rational homology sphere,
then in \cite{Men} we prove the following facts. Using the combinatorics
of the plumbing graph of $M$, one can recover the knot $K_z$ in it.
Then, by the above proposition, from the pair $(M,K_z)$ one can recover
the Alexander polynomial $\Delta_{S^3}(K_f)$ of $f$. It is well-known that
this is equivalent with the equisingular type of the plane singularity $f$.
Moreover, analyzing again the graph of $(M,K_z)$,
one can recover the integer $n$ as well. In particular, from $M$,  we
can recover  not only the geometric genus of $\{f+z^n=0\}$ (fact which is
proved in this article), but also its {\em multiplicity}, and in fact, any
numerical invariant which can be computed from the Newton (or Puiseux) pairs
of $f$ and from the integer $n$ (e.g. even all the equivariant Hodge numbers
associated with the vanishing cohomology of the hypersurface
singularity $g=f+z^n$, or even the
{\em embedded } topological type $(S^5,M)$ of $g$ with its
integral Seifert matrix).

(2) \ref{sum} suggests the following question.
Let $N$ be an integral homology sphere, $L\subset N$ a knot
in it such that $(N,L)$ can be represented by a (negative definite)
plumbing.  Let $(M,K)$ be the $n$-cyclic cover
of $(N,L)$ (branched along $L$) such that $M$ is a rational homology sphere
with $H_1(M)=H$.
Then is it true that $\Delta^H_M(K)(t)=\Delta_N(L)(t)$ ?

The answer is negative: one can construct easily examples (satisfying even the
algebraicity condition) when the identity \ref{sum} fails.
For example, consider $(N,L)$ given by the following splice, respectively
plumbing diagram:

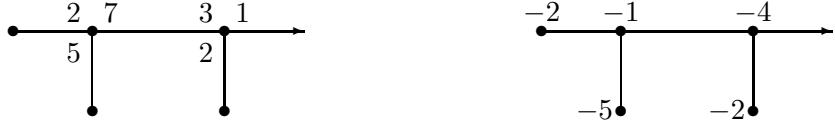
\begin{picture}(400,60)(-30,0)
\put(20,40){\circle*{4}}
\put(50,40){\circle*{4}}
\put(100,40){\circle*{4}}
\put(50,10){\circle*{4}}
\put(100,10){\circle*{4}}
\put(100,40){\vector(1,0){30}}
\put(50,40){\line(0,-1){30}}
\put(100,40){\line(0,-1){30}}
\put(20,40){\line(1,0){80}}

\put(43,47){\makebox(0,0){$2$}}
\put(57,47){\makebox(0,0){$7$}}
\put(93,47){\makebox(0,0){$3$}}
\put(107,47){\makebox(0,0){$1$}}
\put(43,32){\makebox(0,0){$5$}}
\put(93,32){\makebox(0,0){$2$}}

\put(220,40){\circle*{4}}
\put(250,40){\circle*{4}}
\put(300,40){\circle*{4}}
\put(250,10){\circle*{4}}
\put(300,10){\circle*{4}}
\put(300,40){\vector(1,0){30}}
\put(250,40){\line(0,-1){30}}
\put(220,40){\line(1,0){80}}
\put(300,40){\line(0,-1){30}}

\put(220,47){\makebox(0,0){$-2$}}
\put(250,47){\makebox(0,0){$-1$}}
\put(300,47){\makebox(0,0){$-4$}}
\put(240,10){\makebox(0,0){$-5$}}
\put(290,10){\makebox(0,0){$-2$}}
\end{picture}

\noindent Then one can show that e.g. for $n=2$ the identity
$\Delta^H_M(K)(t)=\Delta_N(L)(t)$ fails.

This example also shows that
in the Algebraic Lemma \ref{a3} the assumption $a\geq \deg\Delta$ is crucial.
Indeed, in this example $H=\Z_3$; and in order to determine $\Delta^H_M(K)(t)$,
one needs to compute a sum like in \ref{a3}(a) with
$a=3$, $A=1$ and $\Delta=t^4-t^3+t^2-t+1$ (i.e. with $a<\deg\Delta$).
But for these data, the identity in \ref{a3}(a)  fails.

\subsection{The Reidemeister-Turaev sign-refined torsion.}\label{rt}
Now we will start to compute $\et_{M,\sigma_{can}}(1)$ associated with
$M=M_{(s)}$ and the canonical $spin^c$-structure $\sigma_{can}$ of $M$.
Similarly as above, we write $H=H_1(M)$. Using \ref{sw}(14),
$\et_{M,\sigma_{can}}(1)$ can be determined by the Fourier inversion formula
from $\{\hat{\et}_{M,\sigma_{can}}(\chi)\, :\, \chi\in \hat{H}\setminus \{1\}
\}$. On the other hand, each $\hat{\et}_{M,\sigma_{can}}(\bar{\chi})$ is
given by the limit $\lim_{t\to 1} \hat{P}_{M,\chi,u}(t)$ for some convenient
$u$, cf. \ref{sw}(16).

In the next discussion, the following terminology is helpful.
Fix an integer $1\leq k\leq s$ and a vertex $v'(I):=v'(i_{k+1},\ldots,i_s)
\in \pi^{-1}(v_k)$. Consider the graph $\Gamma\setminus \{v'(I)\}$.
If $\tilde{h}_k=1$ then it has $h_k+2$ connected components: $h_k$ (isomorphic)
subgraphs $\Gamma^{i_k}_-(v'(I))$ ($1\leq i_k\leq h_k$) which contain
vertices at level $k-1$, a string $\Gamma_{st}(v'(I))$
containing a vertex above $\bar{v}_k$, and the component
$\Gamma_+(v'(I))$ which supports the arrow $\{z=0\}$.
Similarly, if $h_k=1$, then $\Gamma\setminus \{v'(I)\}$ has $\tilde{h}_k+2$
connected components, $\Gamma_-(v'(I))$ contains vertices at level $k-1$,
$\Gamma_+(v'(I))$ supports the arrow $\{z=0\}$, and $\tilde{h}_k$
other (isomorphic)
components $\Gamma_{st}^{j_k}(v'(I))$ ($1\leq j_k\leq \tilde{h}_k$),
which are strings, and each of them contains exactly one vertex staying above
$\bar{v}_k$.

Similarly as for the Seifert manifold $\Sigma(p,n,n)$
(see the discussion in \ref{et} after the
diagram), for a large number of characters $\chi$, the limit $\lim_{t\to 1}
\hat{P}_{M,\chi,u}(t)$ is zero. Analyzing the structure of the graph of $M$
and the supports of the characters, one can deduce that a {\em non-trivial }
character $\chi$, with the above limit nonzero, should satisfy one of the
 following  structure properties.

\vspace{2mm}

\noindent {\bf E(asy) case:} The character $\chi$ is born at level $k$
(for some $1\leq k\leq s$) with $\tilde{h}_k>1$. For any vertex $v'\in
\pi^{-1}(v_k)$ one has $\chi(g_{v'})=1$, but there is exactly one vertex
$v'(I):=v'(i_{k+1},\ldots , i_s)\in \pi^{-1}(v_k)$ which is adjacent with the
support of $\chi$. Moreover, $\chi$ is supported by exactly two
components of type $\Gamma_{st}^{j_k}(v'(I))$, say for indices $j_k'$ and
$j_k''$. Let $v'(j_k')$ be the unique vertex in $\Gamma_{st}^{j_k}(v'(I))
\cap \pi^{-1}(\bar{v}_k)$ (similarly for $j_k''$). Then $v'(j_k')$ and
$v'(j_k'')$ are the only vertices $v'$ of the graph of $M$ with
$\delta_{v'}\not=2$ and $\chi(g_{v'})\not=1$. Moreover, $\chi(g_{v'(j_k')})=
\bar{\chi}(g_{v'(j_k'')})=\eta\in \Z_{p_k}^*$.
Therefore, with fixed $(i_{k+1},\ldots, i_s)$ and $(j_k',j_k'')$, there
are exactly $p_k-1$ such characters.

\noindent {\bf D(ifficult) case:}  The character $\chi$ is born at level $k$
(for some $1\leq k\leq s$) with $h_k>1$. For any vertex $v'\in
\pi^{-1}(v_k)$, $\chi(g_{v'})=1$, but there is exactly one vertex
$v'(I):=v'(i_{k+1},\ldots , i_s)\in \pi^{-1}(v_k)$ which is adjacent with the
support of $\chi$. The character $\chi$ is supported by exactly two components
of type $\Gamma_-^{i_k}(v'(I))$, say for indices $i_k'$ and $i_k''$.
Using the previous notations, this means that $\xi_{i_k\,i_{k+1}\cdots i_s}=1$
excepting for $i_k=i_k'$ or $i_k=i_k''$.
(Evidently, $\xi_{i_t\cdots i_s}=1$ for any $t>k$.)
For $t<k$, the values $\xi_{i_t\cdots i_{k+1}\cdots i_s}$ are
un-obstructed. In particular, with indices $(i_{k+1},\ldots, i_s)$ and
$(i_k',i_k'')$ fixed, there are exactly $(a_k-1)\cdot|H_1(M_{(k-1)})|^2$
such characters. Here, $a_k-1$ stands for $ \xi_{i_k'\, i_{k+1}\cdots i_s}=
\bar{\xi}_{i_k''\, i_{k+1}\cdots i_s}\in \Z_{a_k}^*$, and $|H_1(M_{(k-1)})|^2$
for the un-obstructed characters born at level $<k$ on the two branches
corresponding  to $(i_k',i_{k+1},\ldots,i_s)$ and $(i_k'',i_{k+1},\ldots,i_s)$.

\vspace{2mm}

In both cases (E) or (D), if such a character $\chi$ is born at level $k$
(i.e. if satisfies the above characterization for $k$), then we write
$\chi\in B_k$.

\vspace{2mm}

Now, we fix a nontrivial character $\chi$.
Let $S(\chi)$ be the support of $\chi$ and $\bar{S}(\chi)$ its complement.
Then
$$\frac{1}{|H|}\cdot \hat{\et}_{M,\sigma_{can}}(\bar{\chi})=
Loc(\bar{\chi})\cdot Reg(\bar{\chi}),$$
where
$$Loc(\bar{\chi}):=\prod_{v'\in S(\chi)}\, (\chi(g_{v'})-1)^{\delta_{v'}-2}\
\ \ \mbox{and}\ \ \
Reg(\bar{\chi}):=\frac{1}{|H|}\cdot\lim_{t\to 1}\,
\prod_{v'\in \bar{S}(\chi)}\, (t^{w_{v'}(u)}-1)^{\delta_{v'}-2}.$$
We will call $Loc(\bar{\chi})$ the ``local contribution'', while
$Reg(\bar{\chi})$ the ``regularization contribution''.

By the above discussion, $Reg(\bar{\chi})=0$ unless $\chi$ is not of the type
(E) or (D) described above. If $\chi$ is of type (E) or (D)
described as above, then in $\hat{P}_{M,\chi,u}(t)$ (cf. \ref{sw}(15))
one can take $v'(I)$. Moreover, if $\chi\in B_k$,
then by the symmetry of the plumbing graph of $M$, $Reg(\bar{\chi})$
does not depend on the particular choice of $\chi$, but only on the
integer $k$. We write $Reg(k)$ for $Reg(\chi)$ for some (any) $\chi\in B_k$.

In particular,
\begin{equation*}
\et_{M,\sigma_{can}}(1)=\sum_{k=1}^s
Reg(k)\cdot\sum _{\chi\in B_k}\, Loc(\bar{\chi}).
\tag{$\et$}
\end{equation*}

\subsection{Proposition.}\label{Loc} {\em For
any fixed $1\leq k\leq s$ one has:

(E) \ If $h_k=1$ then
$$\sum_{\chi\in B_k}\, Loc(\bar{\chi})=d_k\cdot \frac{
\tilde{h}_k(\tilde{h}_k-1)}{2}\cdot \frac{p_k^2-1}{12}.$$

(D) \ If $\tilde{h}_k=1$ then }
$$\sum_{\chi\in B_k}\, Loc(\bar{\chi})=d_k\cdot \frac{
h_k(h_k-1)}{2}\cdot |H_1(M_{(k-1)})|^2\cdot
\Big[ \,  \frac{a_k^2-1}{12} +(\Delta(f_{(k-1)})^\no)''(1)\, \Big].$$

\begin{proof} (E) $d_k=h_{k+1}\cdots h_s$ is the cardinality of the
index set $(i_{k+1},\ldots,,i_s)$, $\tilde{h}_k(\tilde{h}_k-1)/2$ is the
number of possibilities to choose the indices $(j_k',j_k'')$.
The last term comes from a formula of type \ref{g8}($**$), where the sum
is over $\eta\in \Z_{p_k}^*$.

In the case (D), $d_kh_k(h_k-1)/2$ has the same interpretation.
Fix the branch $(i_k',i_{k+1},\ldots,i_s)$ and consider the  sum over all the
characters born at level $<k$. Then \ref{sum},  applied for $(M_{(k-1)},
K_{(k-1)})$ as a covering of $(S^3,f_{(k-1)}=0)$, provides
$|H_1(M_{(k-1)})|^2\cdot \Delta(f_{(k-1)})(t)/(t-1)$
evaluated at $t=\xi_{i_k'\, i_{k+1}\cdots i_s}$. The same is true for
the other index $i_k''$. Then apply \ref{a4}(3) for $\Delta=\Delta(f_{(k-1)})$
and $a=a_k$. This can be applied because of \ref{a2} and \ref{a1}(6).
\end{proof}

\subsection{The ``regularization contribution'' $Reg(k)$.}\label{refl} \
Fix a character $\chi\in B_k$ of type (E) or (D) as in \ref{rt}.
Recall that one can
take $u=v'(I)$. Consider the connected components of $\Gamma\setminus \{v'(I)\}
$ (as in \ref{rt}), where we add to each component an arrow corresponding
to the edge which connects the component to $v'(I)$. For these graphs, if
one applies \ref{s0e}(12), one gets that $Reg(k)$ is

\[
Reg(k)=\left\{
\begin{array}{cc}{-\det(\Gamma_-)\cdot \det(\Gamma_{st})^{\tilde{h}_k-2}
\cdot\det(\Gamma_+)/\det(\Gamma)}&\mbox{case (E)}\\
& \\
-\det(\Gamma_-)^{h_k-2}\cdot \det(\Gamma_{st})
\cdot\det(\Gamma_+)/\det(\Gamma)&\mbox{case (D)}.
\end{array}
\right.
\]
Let $I$ be the intersection matrix of $M$, and
$I^{-1}_k:=I^{-1}_{v'v'}$ for any $v'\in \pi^{-1}(v_k)$.
Then, by the formula which provides the entries of an inverse matrix,
one gets that
\[
\det \Gamma\cdot I_{v'v'}^{-1}=
\left\{
\begin{array}{cc}\det(\Gamma_-)
\det(\Gamma_{st})^{\tilde{h}_k}\det(\Gamma_+)& \mbox{case (E)}\\
&\\
\det(\Gamma_-)^{h_k}\det(\Gamma_{st})\det(\Gamma_+)&\mbox{case (D)}.
\end{array}
\right.
\]
These two facts combined show that
\[
Reg(k)=
\left\{
\begin{array}{lll}
-I^{-1}_{k}/\det(\Gamma_{st})^2 &=
-I_k^{-1}/p_k^2& \mbox{case (E)}\\
& & \\
-I^{-1}_{k}/\det(\Gamma_{-})^2&=
-I_k^{-1}/(a_k\cdot |H_1(M_{(k-1)})|)^2 & \mbox{case (D)}
\end{array}
\right.
\]
since in case (E) $|\det(\Gamma_{st})|=p_k$ e.g. from \ref{g2}(a)), and
in case (D) $|\det(\Gamma_{-})|=a_k\cdot |H_1(M_{(k-1)})|$ by \ref{detg}.

This, together with \ref{rt} ($\et$) and \ref{Loc} transform into
\begin{equation*}
\et_{M,\sigma_{can}}(1)=-\sum_{k=1}^s I^{-1}_k\cdot d_k\cdot A_k/2,
\tag{$*$}
\end{equation*}
where $A_k$ is defined in \ref{pr}(k) in terms of the numerical invariants of
$f_{(k)}$.

\subsection{The computation of $I^{-1}_k$.}\label{ii} For any $l\geq k$, let
$I^{-1}_k(M_{(l)})$ be the $(v',v')$-entry of the inverse of the intersection
form $I(M_{(l)})$ associated with $M_{(l)}$, where $v'$ is any vertex
above $v_k$. E.g. $I^{-1}_k(M_{(s)})$ is $I^{-1}_k$ used above.

By \ref{s0d}(5), $-I^{-1}_k(M_{(l)})=Lk_{M_{(l)}}(g_{v'},g_{v'})$.
If $l=k$, by \ref{s5}(12) this is $Lk_{\Sigma(p_k,a_k,n/d_k)}(O,O)$,
hence by \ref{g3} (L\#) it is $np_ka_k/(d_kh_k^2\tilde{h}_k^2)$.

Next, assume that $l>k$. If $h_l=1$, then the splicing
\[
M_{(l)}=h_lM_{(l-1)}\coprod \Sigma(p_l,a_l,n/d_l)
\]
 is trivial (with $o_1=o_2=1$ and $k_1=k_2=0$), hence by \ref{s5}(10) one gets
 \[
 -I^{-1}_k(M_{(l)})=-I^{-1}_k(M_{(l-1)}).
 \] If $h_l>1$, then \ref{g4} and
an iterated application of \ref{s5}(10) and  \ref{g3} (L\#)  give
$$-I^{-1}_k(M_{(l)})=-I^{-1}_k(M_{(l-1)})-\Big(\frac{a_kp_k}{\tilde{h}_kh_k}
\cdot \frac{p_{k+1}}{h_{k+1}}\cdots \frac{p_{l-1}}{h_{l-1}}\Big)^2\cdot
\frac{np_l(h_l-1)}{d_la_lh_l^2}.$$
Indeed, by \ref{s0d}(5) and \ref{s5}
$$-I^{-1}_k(M_{(l)})=-I^{-1}_k(M_{(l-1)})-\Big(Lk_{M_{(l-1)}}(g_{v'},K_{(l-1)})
\Big)^2\cdot \frac{np_l(h_l-1)}{d_la_lh_l^2},$$
and, again by \ref{s5},  $Lk_{M_{(l-1)}}(g_{v'},K_{(l-1)})$ equals
$$Lk_{\Sigma(p_k,a_kn/d_k)}(O,Z)\cdot
Lk_{\Sigma(p_{k+1},a_{k+1}n/d_{k+1})}(K_2^{(i)},Z)\cdots
Lk_{\Sigma(p_{l-1},a_{l-1}n/d_{l-1})}(K_2^{(i)},Z).$$
%By induction, the above formula can be written as
%$$-I^{-1}_k(M_{(s)})=\frac{na_kp_k}{d_k\tilde{h}_kh_k}\, \Big[
%1-\frac{a_kp_kp_{k+1}(h_{k+1}-1)}{a_{k+1}h_{k+1}}-
%\frac{a_kp_kp_{k+1}^2p_{k+2}(h_{k+2}-1)}{a_{k+2}h_{k+1}h_{k+2}}-\cdots $$
%$$\hspace{2cm}\cdots -
%\frac{a_kp_kp_{k+1}^2\cdots p_{s-1}^2p_s(h_s-1)}{a_sh_{k+1}h_{k+2}\cdots h_s}
%\Big].$$

\subsection{The splicing formula for $\et_{M\sigma_{can}}(1)$.}\label{splet} \
Using \ref{refl}($*$) and \ref{ii} (and $d_s=1$),  we can write
$$\et_{M_{(s)},\sigma_{can}}(1)-h_s\cdot
\et_{M_{(s-1)},\sigma_{can}}(1)=\frac{na_sp_s}{2\tilde{h}_s^2h_s^2}\cdot A_s$$
$$+\sum_{k=1}^{s-1}(-I_k^{-1}(M_{(s)}))\cdot h_{k+1}\cdots h_s\cdot A_k/2-
h_s\sum_{k=1}^{s-1}(-I_k^{-1}(M_{(s-1)}))\cdot h_{k+1}\cdots h_{s-1}
\cdot A_k/2$$
$$=\frac{na_sp_s}{2\tilde{h}_s^2h_s^2}\cdot A_s
-\sum_{k=1}^{s-1} h_{k+1}\cdots h_s\cdot \frac{a_k^2p_k^2\cdots p_{s-1}^2}
{\tilde{h}_k^2h_k^2\cdots h_{s-1}^2}\cdot \frac{np_s(h_s-1)}{a_sh_s^2}\cdot
A_k/2.$$
But, by \ref{g2}(e)  one has
$$\frac{na_sp_s}{2\tilde{h}_s^2h_s^2}\cdot A_s=
\et_{\Sigma(p_s,a_s,n),\sigma_{can}}(1)+
\frac{na_sp_s}{2\tilde{h}_s^2h_s^2}\cdot \frac{h_s(h_s-1)}{a_s^2}\cdot
(\Delta(f_{(s-1)})^\no)''(1).$$
Therefore, (using also $(h_s-1)/\tilde{h}^2=h_s-1$) one gets
$${\mathcal O}(\et_{\cdot,\sigma_{can}}(1))=
\et_{M_{(s)},\sigma_{can}}(1)-h_s\cdot \et_{M_{(s-1)},\sigma_{can}}(1)
-\et_{\Sigma(p_s,a_s,n),\sigma_{can}}(1)$$
$$=\frac{np_s(h_s-1)}{2h_sa_s}\Big[
(\Delta(f_{(s-1)})^\no)''(1)-\sum_{k=1}^{s-1}\, \frac{a_k^2p_k^2\cdots
p_{s-1}^2}{\tilde{h}_k^2h_k^2h_{k+1}\cdots h_{s-1}}\cdot A_k\Big].$$
Notice that this inductive formula (together with \ref{g2}(e))
shows that the numerical function $n\mapsto \et(n)$ (where $\et(n)$ denotes
the sign-refined Reidemeister-Turaev torsion of the link of $f+z^n$,
associated with $\sigma_{can}$) can be written in the form
$n\mapsto np_1(n)+p_2(n)$, where $p_i(n)$ are periodic functions.
Similar result for the signature of the Milnor fiber was obtained by
Neumann \cite{Ncy}

\vspace{2mm}

Now, using \ref{g9} for $M=M_{(s)}$, and \ref{pr}(k) for $l=s-1$,  one has the
following  consequences:

\subsection{Theorem.}\label{swadd} {\em The additivity obstruction
${\mathcal O}(\ssw^0_\cdot  (\sigma_{can}))=0$, in other words:
$$\ssw^0_{M_{(s)}}(\sigma_{can})=h_s\cdot
\ssw^0_{M_{(s-1)}}(\sigma_{can})+
\ssw^0_{\Sigma(p_s,a_s,n)}(\sigma_{can}).$$
In particular, by induction for $M=M_{(s)}$ one gets:}
$$\ssw^0_{M}(\sigma_{can})=\sum_{k=1}^s\, d_k\cdot
\ssw^0_{\Sigma(p_k,a_k,n/d_k)}(\sigma_{can}).$$

\subsection{Corollary.}\label{conj} {\em Consider the
hypersurface singularity $g(x,y,z)=f(x,y)+z^n$, where $f$ is an
irreducible plane curve singularity.
Assume that its link $M$ is
a rational homology sphere. If $\sigma(g)$ denotes the signature
of the Milnor fiber of $g$, then
$$-\ssw^0_M(\sigma_{can})=\sigma(g)/8.$$
In particular, the geometric genus of $\{g=0\}$ is topological and it is
given by
$$p_g=\ssw^0_M(\sigma_{can})-(K^2+\#{\mathcal V})/8,$$
where the invariant $K^2+\#{\mathcal V}$ (associated with any plumbing
graph of $M$) is defined in the introduction.}

\begin{proof} By \ref{swadd} an \cite{NDedI}(3.2) (cf. also with
\ref{pr}(e)), we only have to show that
$$-\ssw^0_{\Sigma(p_k,a_k,n/d_k)}(\sigma_{can})=
\sigma(p_k,a_k,n/d_k)/8$$ for any $k$. But this follows from \cite{NN}
(section 6). \end{proof}

For an explicit formula of
$\ssw^0_{\Sigma(p_k,a_k,n/d_k)}(\sigma_{can})$, see \cite{NN} or \cite{NNII}.

\end{document}